\documentclass{article}
\usepackage[utf8]{inputenc}

\usepackage[margin=3cm]{geometry}
\usepackage{amsfonts}
\usepackage{amsthm}
\usepackage{amsmath}
\usepackage{mathtools}
\usepackage{amssymb}
\usepackage{multicol}
\usepackage{cite}
\usepackage{caption}
\usepackage{subcaption}
\usepackage{graphicx}
\usepackage{tcolorbox}
\usepackage{dsfont}
\usepackage{xcolor}
\usepackage{tikz-cd}
\usepackage{subfiles}
\usepackage{longtable}
\usepackage{hyperref}
\usepackage{bbm}
\newtheorem{theorem}{Theorem}
\newtheorem*{theorem*}{Theorem}
\newtheorem{corollary}[theorem]{Corollary}

\newtheorem{lemma}[theorem]{Lemma}
\newtheorem{defn}[theorem]{Definition}
\newtheorem{propn}[theorem]{Proposition}
\newtheorem*{propn-non}{Proposition}

\makeatletter
\def\blfootnote{\gdef\@thefnmark{}\@footnotetext}
\makeatother
\newcommand{\B}[1]{\mathbb{#1}}
\newcommand{\C}[1]{\mathcal{#1}}
\newcommand{\oop}{\preccurlyeq}
\newcommand{\is}{\mathbbm{1}}

\DeclareMathOperator{\erfc}{erfc}

\usepackage[skip=7.5pt plus1pt, indent=20pt]{parskip}

\title{Critical Drift for Brownian Bees and a Reflected Brownian Motion Invariance Principle}
\author{Jacob Mercer\footnote{\texttt{jacob.mercer@maths.ox.ac.uk}, Department of Statistics, University of Oxford}}
\date{}

\begin{document}
\maketitle
\begin{abstract}
    $N$-Brownian bees is a branching-selection particle system in $\B{R}^d$ in which $N$ particles behave as independent binary branching Brownian motions, and where at each branching event, we remove the particle furthest from the origin. We study a variant in which $d=1$ and particles have an additional drift $\mu\in\B{R}$. We show that there is a critical value, $\mu_c^N$, and three distinct regimes (sub-critical, critical, and super-critical) and we describe the behaviour of the system in each case. In the sub-critical regime, the system is positive Harris recurrent and has an invariant distribution; in the super-critical regime, the system is transient; and in the critical case, after rescaling, the system behaves like a single reflected Brownian motion.

    We also show that the critical drift $\mu_c^N$ is in fact the speed of the well-studied $N$-BBM process, and give a rigorous proof for the speed of $N$-BBM, which was missing in the literature. 
\end{abstract}
\section*{Introduction} \blfootnote{This publication is based on work supported by the EPSRC Centre for Doctoral Training in Mathematics of Random Systems: Analysis, Modelling, and Simulation (EP/S023925/1)}

In the early 2000's, the physicists Brunet and Derrida pioneered the study of branching particle systems with selection. In these so called \textit{Brunet-Derrida} or \textit{Branching-Selection} systems, we have a number of particles which can branch and have offspring, and also some \textit{selection mechanism/rule} which removes individuals from the population, keeping the population size constant. The particles removed at branching events are typically the particle of least `fitness' according to some fitness function; for example $f(x)=-\|x\|$, in which case it is always the particle furthest from the origin which is removed. We can also remove particles from the system by measuring fitness with respect to the whole system; for example, removing the particle furthest from the barycentre (as in \cite{bbb}). The selection mechanism therefore introduces a `survival-of-the-fittest' phenomenon into the model. 

A well known branching-selection particle system is the $N$-branching Brownian motion ($N$-BBM) which is studied in \cite{selecJBOT}, \cite{hydroNBBM}, \cite{maillard}. In the $N$-BBM, $N$ particles behave as independent Brownian motions on the real line, each branching into 2 particles at rate 1, and at each branching event, the leftmost particle is removed in order to keep a fixed-size population of $N$ particles. This mechanism  introduces a \textit{selective pressure} which `pushes' the cloud of particles to the right. Writing $Z(t)=(Z_1(t),...,Z_N(t))\in\B{R}^N$ for the positions of the particles at time $t$, then we say that $Z$ has velocity $v_N$ if the following limits exist and are equal:
$$
v_N = \lim_{t\to\infty}\min_{1\leq i\leq N} \frac{Z_i(t)}{t} = \lim_{t\to\infty}\max_{1\leq i\leq N} \frac{Z_i(t)}{t}
$$
It is well known in the literature (see for example Section 1.3 of \cite{maillard} or Remark 1.6 of \cite{shape}) that the velocity of the $N$-BBM is given by:
$v_N = \sqrt{2} - \frac{\pi^2}{\sqrt{2}(\log N)^2} + o\left(\frac{1}{(\log N)^2}\right).$
This is very similar to B\'{e}rard and Gou\'{e}r\'{e} (\hspace{1sp}\cite{berardGouere}), who prove an analogous result for the $N$-branching random walk. A complete proof for the case of the $N$-BBM is given in the appendix.

Another branching-selection particle system which has attracted attention lately is called `Brownian bees', and is studied by Berestycki, Brunet, Nolen, and Pennington in \cite{afree}, \cite{beesBBNP}. In the $N$-Brownian bees particle system, $N$ particles behave as Brownian motions, each branching at rate $1$. At each branching event, we delete the particle which is furthest from the origin. The $N$-Brownian bees system is defined in $\B{R}^d$ for any number of dimensions, $d$, but in this paper we specifically study $d=1$. 

A number of things are known about the $N$-Brownian bees system. Firstly, the system is positive Harris recurrent (Proposition 6.5 in \cite{beesBBNP}). It is also known that in the large population limit, the empirical distribution of particles converges to the solution of a free boundary problem. The aim of this paper is to study the $N$-Brownian bees system in which each particle also has a drift of $\mu$. 

Intuitively, when the drift is small, the system shouldn't behave much differently to the $N$-Brownian bees system without drift. Essentially, the selection rule which deletes the furthest particle from the origin outweighs the effect of the drift and succeeds in keeping all particles near to the origin. Alternatively, when the drift is sufficiently large, it will win against the selection mechanism and particles will be drifting off to $\pm\infty$. 

We will therefore prove that there is a critical drift $\mu_c^N$ such that when $|\mu| < \mu_c^N$, the $N$-Brownian bee system is positive Harris recurrent, and when $|\mu| > \mu_c^N$, the system is transient in the sense that each particle diverges to $\frac{\mu}{|\mu|}\infty$ as $t\to\infty$. 

Now we can observe that when all $N$ particles of a Brownian bees system lie in $(-\infty,0]$, the system has the same behaviour as $N$-BBM with added drift $\mu$, which has asymptotic velocity $\mu + v_N$. Therefore we should expect the critical drift $\mu_c^N$ to be:
$$\mu_c^N = v_N = \sqrt{2} - \frac{\pi^2}{\sqrt{2}(\log N)^2} + o\left(\frac{1}{(\log N)^2}\right)$$

In the critical case where $|\mu|=\mu_c^N$, the system becomes null-recurrent and, under a suitable rescaling, the $N$ particles act like a single point mass behaving as a reflected Brownian motion. Our main result is  the following theorem:

\begin{theorem} \label{mainThem} Let $X^{N,\mu}=(X^{N,\mu}(t))_{t\geq 0}$ be a one dimensional $N$-Brownian bee system with drift $\mu\in \B{R}$, and let $\mu_c^N=v_N$. Then:
\begin{enumerate}
    \item \underline{Sub-Critical Case}: when $|\mu| < \mu_c^N$, then for any $t_0>0$, the process $(X^{N,\mu}(nt_0))_{n\in \B{N}}$ is positive Harris recurrent and therefore has a unique stationary probability measure $\pi^{N,\mu}$ on $\B{R}^N$ so that for any $\chi\in \B{R}^N$
    \begin{equation*}
        \lim_{t\to\infty}\sup_{C}|\B{P}_{\chi}(X^{N,\mu}(t)\in C)-\pi^{N,\mu}(C)|_{TV}=0.
    \end{equation*}    
    where $|\cdot|_{TV}$ represents the total variation norm, the supremum is taken over all measurable sets, and $\B{P}_\chi(\cdot):=\B{P}(\cdot|X^{N,\mu}(0)=\chi)$.
    \item \underline{Critical Case}: When $|\mu |=\mu_c^N$, there is an invariance principle
    $$  
    \left(m^{-1/2}X^{N,\mu}(mt)\right)_{t\geq 0} \xrightarrow[m\to\infty]{d} \left(\frac{\mu}{|\mu|}\beta^{-1/2}\sigma |B(t)|\underline{1}\right)_{t\geq 0}.
    $$
    in the Skorokhod topology on $\C{D}([0,\infty),\B{R}^N)$, the space of cadlag functions $[0,\infty)\to \B{R}$, where $B$ is a standard Brownian motion, $\underline{1}=(1,\ldots,1)\in \B{R}^N$, and $\beta,\sigma$ are constants which we will define later. 
    \item \underline{Super-Critical Case}: When $|\mu|>\mu_c^N$, $X^{N,\mu}$ is transient in the sense that each particle $X_i^{N,\mu}$ diverges to $\frac{\mu}{|\mu|}\infty$ as $t\to\infty$, and further
    $$\lim_{t\to\infty}\frac{X^{N,\mu}_1(t)}{t}=\lim_{t\to\infty}\frac{X^{N,\mu}_N(t)}{t} = \left(|\mu| - \mu_c^N\right)\frac{\mu}{|\mu|} \quad \text{a.s.}$$
    where $X_i^{N,\mu}$ denotes the $i$\textsuperscript{th} smallest particle of 
    $X^{N,\mu}$.
\end{enumerate}
\end{theorem}

The subcritical and supercritical cases of Theorem \ref{mainThem} were also independently proven in the master's thesis of Flynn (\hspace{1sp}\cite{flynn}).

As we noted above, the behaviour of the $N$-Brownian bees system without drift is known when the number of particles tends to $\infty$. In particular for $d=1$, it is known that as $N\to\infty$, the empirical distribution of particles converges to the solution of a free boundary problem given by
\begin{align} \label{beesFBP1}
\begin{cases}
    u_t=\frac{1}{2}u_{xx} + u  & x\in(-R_t,R_t), t>0, \\
    \int_{-R_t}^{R_t} u(s,t)ds=1 & t>0, \\
    u(x,t)=0 &  x\notin (R_t,R_t), t>0, \\
    u \text{ continuous}.
\end{cases}
\end{align}
with an initial condition based on the initial distribution of particles. This result is proven by Berestycki et al. in \cite{afree}, \cite{beesBBNP}. We conjecture here that for Brownian bees with sub-critical drift $\mu$, the empirical distribution converges to the solution of the free boundary problem: 
\begin{align}
\begin{cases}
    u_t=\frac{1}{2}u_{xx} - \mu u_x + u  & x\in(-R_t,R_t), t>0, \\
    \int_{-R_t}^{R_t} u(s,t)ds=1 & t>0, \\
    u(x,t)=0 &  x\notin (R_t,R_t), t>0, \\
    u \text{ continuous}.
\end{cases}
\end{align}

The methods of proof used by Berestycki et al. rely on the symmetry of the system, and view it essentially as a system of reflected Brownian motions on $[0,\infty)$. However these methods no longer work when drift is introduced, since $|B_t + \mu t|$ does not behave like a reflected Brownian motion and in fact is not even Markovian. Therefore this conjecture remains unproven.

\section{Some useful couplings of $N$-BBM and Brownian bees}
Results about the $N$-BBM are often easier to prove than results about the $N$-Brownian bees system due to the fact that the $N$-BBM is translation invariant. A number of our proofs will use properties of the $N$-BBM and then construct a coupling between $N$-BBM and $N$-Brownian bees in order to yield the desired results. Therefore in the first section, we introduce a number of couplings which we will need in the proof of Theorem \ref{mainThem}.

We begin by defining the $N$-BBM and $N$-Brownian bees processes exactly. In different proofs, it will make sense to either describe the process by $N$ particles with intrinsic labels which do not change, or to label the particles by order, so we give both definitions here. For an $\B{R}^N$-process, define $\Theta^N:\B{R}^N \to \B{R}^N$ to be the function which ranks the components of the vector, breaking ties arbitrarily. So the $N$-particle processes $(A(t))_{t\geq 0}$ and $(\Theta^N(A(t))_{t\geq 0}$ have the same empirical measure but different labellings of components. 

\vspace*{1em}

\begin{defn} \label{nbbmDef} (\textbf{$N$-BBM}) Let $(W^N(t))_{t\geq 0}=\left((W^N_1(t))_{t\geq 0},\ldots,(W^N_N(t))_{t\geq 0}\right)$ be the trajectories of $N$ particles on the real line. The particles behave as $N$ independent Brownian motions, and at rate $N$ the leftmost particle in the system jumps to the position of a uniformly randomly chosen particle. At any time $t\geq 0$, let $Z^N(t):=\Theta^N(W^N(t))$; so $Z^N$ is the process with ordered components. In the rest of this paper, we will call an $N$-BBM any cadlag stochastic process which  has the same empirical distribution as $W^N$ or $Z^N$
\end{defn}

\vspace*{1em}

\begin{defn} \label{beesDef} (\textbf{$N$-Brownian Bees}) Let $(V^N(t))_{t\geq 0} = \left((V^N_1(t))_{t\geq 0},\ldots,(V^N_N(t))_{t\geq 0}\right)$ be the trajectories of $N$ particles on the real line. The particles behave as $N$ independent Brownian motions, and at rate $N$ the particle furthest from the origin jumps to the position of a uniformly randomly chosen particle. At any time $t\geq 0$, let $X^N(t):=\Theta^N(V^N(t))$; so $X^N$ is the process with ordered components. In the rest of this paper, we will call an $N$-Brownian bees any cadlag stochastic process which has the same empirical distribution as $V^N$ or $X^N$.
\end{defn}

\paragraph{Remark:} We could equivalently think of the $N$-BBM (resp. $N$-Brownian bees) as being a process in which  each particle branches at unit rate and at each branching event the leftmost particle (resp. furthest particle from the origin) is deleted, giving its label to the newly born particle. 

We also define a notion of `ordering' for particle systems.

\vspace*{1em}

\begin{defn} Let $A=(a_1,\ldots,a_m)\in \B{R}^m$ and $B=(b_1,\ldots, b_n)\in \B{R}^n$. Then $A\oop B$ (we say that ``$A$ lies to the left of $B$'', or ``$A$ is less than $B$ in the sense of $\oop$'') if and only if
$$|\{i:a_i \geq c\}| \leq |\{i: b_i \geq c\}| \; \forall c\in \B{R}.$$
When $m=n$, this means that the empirical cumulative distribution functions are pointwise ordered.
\end{defn}

Now our first result will prove that both the $N$-BBM and the $N$-Brownian bees can be constructed from i.i.d. Brownian motions, a Poisson point process, and i.i.d. uniform random variables on $\{1,2,\ldots,N\}$. This representation of $N$-BBM and $N$-Brownian bees will then make proving couplings with certain properties very straightforward. The first natural coupling this result gives is that the $N$-BBM process is monotonic as a function of its initial configuration, and the second is that an $N$-Brownian bees process can be coupled to always stay to the left of an $N$-BBM process.

\vspace*{1em}

\begin{propn} \label{generalCouple} Let $B=((B^{j}(t))_{t\geq 0})_{1\leq j\leq N}$ be a family of $N$ i.i.d Brownian motions, $\nu$ be an initial configuration of $N$ particles, $Q=(Q(t))_{t\geq 0}$ be a Poisson point process of intensity $N$, and $I=(I_i)_{i\in \B{N}}$ be a sequence of i.i.d random variables which are uniform on $\{1,2,\ldots,N\}$. Then there exist functions $\Upsilon$ and $\Xi$ with values in $\mathbb R^N$ such that:
$$(X^N(t))_{t\geq 0}=\Theta^N\left( \Upsilon\left(\nu,B, Q, I,t\right)_{t\geq 0}\right)$$ is an $N$-Brownian bees process and
$$(Z^N(t))_{t\geq 0}=\Theta^N\left( \Xi\left(\nu,B, Q, I,t\right)_{t\geq 0}\right)$$
is an $N$-BBM process.

Moreover, given configurations $\nu$ and $\nu'$ such that $\nu \oop \nu'$, we have
$$\Upsilon(\nu, B, Q, I,t) \oop \Xi(\nu, B, Q, I,t)\oop \Xi(\nu', B, Q, I,t) \; \forall t\geq 0$$
\end{propn}

\begin{proof}
    We will construct the processes inductively. Let $\tau_0=0$, and recursively, define $$\tau_{i+1}:=\inf\{t> \tau_i:Q(t)\neq Q(t-)\}$$ to be the discontinuities of the Poisson point process $Q$. These times will be the branching times of the processes. For any initial configuration $\eta$, define
    $$\Upsilon(\eta, B,Q,I,0) = \Xi(\eta, B,Q,I,0)=\eta .$$

    Now suppose for induction that $\Upsilon$ and $\Xi$ are defined up to a time $\tau_i$ so that $\Upsilon\left(\nu,B, Q, I,t\right)_{0\leq t\leq \tau_i}$ is an $N$-Brownian bees process and $\Xi\left(\nu,B, Q, I,t\right)_{0\leq t\leq \tau_i}$ and $\Xi\left(\nu',B,Q,I,t\right)_{0\leq t \leq \tau_i}$ are $N$-BBM processes with $$\Upsilon\left(\nu,B, Q, I,t\right)\oop \Xi\left(\nu,B, Q, I,t\right)\oop \Xi\left(\nu',B,Q,I,t\right) \text{ for } t\in [0,\tau_i].$$

    We will then define $\Upsilon$ and $\Xi$ on $[\tau_i,\tau_{i+1}]$. Let $\Upsilon_j(\nu, B,Q,I,\tau_i)$ (resp. $\Xi_j(\nu, B,Q,I,\tau_i)$) denote the $j$\textsuperscript{th} smallest particle of $\Upsilon(\nu, B,Q,I,\tau_i)$ (resp. $\Xi(\nu,B,Q,I,\tau_i)$). Then for $t\in [\tau_i,\tau_{i+1})$, we will drive the particles $\Upsilon_j(\nu,B,Q,I,t)$, $\Xi_j(\nu,B,Q,I,t)$, and $\Xi_j(\nu',B,Q,I,t)$ by the same Brownian motion $(B^{j}(t))_{\tau_i \leq t\leq \tau_{i+1}}$. Since $\Upsilon(\nu,B,Q,I,\tau_i)\oop \Xi(\nu,B,Q,I,\tau_i)\oop \Xi(\nu',B,Q,I,\tau_i)$ implies that $\Upsilon_j(\nu,B,Q,I,\tau_i)\leq \Xi_j(\nu,B,Q,I,\tau_i)\leq \Xi_j(\nu',B,Q,I,\tau_i)$ for all $j\in \{1,2,\ldots,N\}$, then certainly:
    $$\Upsilon_j(\nu,B,Q,I,\tau_i)+B^{j}(t)-B^{j}(\tau_i)\leq \Xi_j(\nu,B,Q,I,\tau_i)+B^{j}(t)-B^{j}(\tau_i)\leq \Xi_j(\nu',B,Q,I,\tau_i)+B^{j}(t)-B^{j}(\tau_i),$$
    for all $j\in \{1,2,\ldots,N\}$, and so we have $\Upsilon(\nu,B,Q,I,t)\oop \Xi(\nu,B,Q,I,t)\oop \Xi(\nu',B,Q,I,t)$ for $t\in [\tau_i,\tau_{i+1})$. 

    Now let us define two operators $k,l:\B{R}^N \times \{1,2,\ldots,N\} \to \B{R}^N$. For a vector $v=(v_1,v_2,\ldots,v_N)$ with ordered components $v_1\leq v_2 \leq \cdots \leq v_N$ we define:
    $$l(v,i)=(v_2,v_3,\ldots,v_{i-1},v_i,v_i,v_{i+1},\ldots,v_N),$$
    and 
    $$k(v,i)=\begin{cases}
        (v_2,v_3,\ldots,v_{i-1},v_i,v_i,v_{i+1},\ldots,v_N) & \text{ if }|v_1|\geq |v_N| \\
        (v_1,v_2,\ldots, v_{i-1},v_i,v_i,v_{i+1},\ldots,v_{N-1}) & \text{ if }|v_N| >|v_1|
    \end{cases}$$
    In words, if $v$ describes the positions of $n$ particles, then $l(v,i)$ duplicates the $i$\textsuperscript{th} smallest particle of $v$ whilst killing the leftmost particle, and $k(v,i)$ duplicates the $i$\textsuperscript{th} smallest particle of $v$ whilst killing the particle of largest magnitude. Note that by definition, we certainly have that $k(v,i)\oop l(v,i)\oop l(v',i)$ for any $v, v' \in \B{R}^N$ such that $v\oop v'$ and $i\in \{1,2,\ldots,N\}$. Finally, we can define $$\Upsilon(\nu,B,Q,I,\tau_{i+1})=k(\Theta^N(\Upsilon(\nu,B,Q,I,\tau_{i+1}-)),I_{i+1}),$$ and $$\Xi(\nu,B,Q,I,\tau_{i+1})=l(\Theta^N(\Xi(\nu,B,Q,I,\tau_{i+1}-)),I_{i+1}),$$ so that it clearly follows that:
    $$\Upsilon(\nu,B,Q,I,\tau_{i+1})\oop \Xi(\nu,B,Q,I,\tau_{i+1})\oop \Xi(\nu',B,Q,I,\tau_{i+1}).$$

    Hence by induction, this ordering holds for all $t\geq 0$. It is obvious from the construction that $\Upsilon(\nu,B,Q,I,t)_{t\geq 0}$ is indeed an $N$-Brownian bees process and $\Xi(\nu,B,Q,I,t)_{t\geq 0}$ is indeed an $N$-BBM process.
\end{proof}

Now we note that if we want to consider an $N$-BBM or an $N$-Brownian bees process with added drift of $\mu$, we can simply consider it as $\Upsilon(\nu,B^\mu,Q,I,t)_{t\geq 0}$ or $\Xi(\nu,B^\mu,Q,I,t)_{t\geq 0}$ respectively, where $B^\mu$ is a family of i.i.d Brownian motions each with drift $\mu$. It immediately follows that the coupling still holds. 

From here on, $Z^{N,\mu}$ will be an $N$-BBM process in which each particle has an additional drift of $\mu\in\B{R}$ on each particle, and similarly $X^{N,\mu}$ will be an $N$-Brownian bees process with an additional drift of $\mu \in \B{R}$ on each particle. 

The next result, Lemma \ref{absCouple}, gives another coupling between an $N$-Brownian bees process and an $N$-BBM. It essentially says that, on time intervals when no particle hits the origin, we can couple the absolute values of $N$-Brownian bees with $N$-BBM. To prove this, we will give an alternative way to represent the $N$-Brownian bees as a function of an initial configuration, a Poisson point process, $N$ i.i.d Brownian motions, and a sequence of uniform random variables. The ordering of this coupling will break down as soon as a particle hits the origin, but this will not be a problem for how we will use the result.

\vspace*{1em}

\begin{lemma} \label{absCouple} Let $\mu \leq 0$, and  let $B$, $\nu$, $Q$, $I$, and $\Xi$ be as in Proposition \ref{generalCouple}. Then there exists a function $\tilde{\Upsilon}$  such that we can write:
$$(X^{N,\mu}(t))_{t\geq 0}=\Theta^N\left(\tilde{\Upsilon}(\nu,B^\mu,Q,I,t)_{t\geq 0}\right).$$

so that, if $T:=\inf\{t\geq 0:X_i^{N,\mu}(t)=0 \text{ for some i}\}$ is the first hitting time of the origin by a particle of $X^{N,\mu}$, and $\tilde{\nu}\oop -|\nu|$ (where $|\nu|$ indicates taking absolute value element-wise), then:
$$\Xi(\tilde{\nu},B^\mu,Q,I,t) \oop -|\tilde{\Upsilon}(\nu,B^\mu,Q,I,t)|\; \text{ for }t\leq T.$$  
\end{lemma}

\begin{proof} As with $\Upsilon(\nu,B^\mu,Q,I,t)$, we will define $\tilde{\Upsilon}(\nu,B^\mu,Q,I,t)_{t\geq 0}$ recursively. Let $\tilde{\Upsilon}(\nu,B^\mu,Q,I,0)=\nu$, and suppose for induction that we have constructed $\tilde{\Upsilon}(\nu,B^\mu,Q,I,t)$ up to time $\tau_i$ so that $\Xi(\tilde{\nu},B^\mu,Q,I,t)\oop -|\tilde{\Upsilon}(\nu,B^\mu,Q,I,t)|$ for $t\in [0,\tau_i]$.  We will now construct the process up to time $\tau_{i+1}$. Similarly to the previous proof, $\tilde{\Upsilon}_j(\nu,B^\mu,Q,I,t)$ denotes the $j$\textsuperscript{th} smallest particle of $\tilde{\Upsilon}(\nu,B^\mu,Q,I,t)$ and define 
\[S_{i,j} := \text{sign} \left( \tilde{\Upsilon}_j(\nu,B^\mu,Q,I,\tau_i) \right) \in\{-1,+1\}.\]

Let us write the Brownian motions with drift $\mu$, $B^\mu$ as $B^{j}(t)+\mu t$ for each $1\leq j\leq N$. Then for $t\in [\tau_i,\tau_{i+1})$, we will drive the particle $\tilde{\Upsilon}_j(\nu,B^\mu,Q,I,\tau_i)$ by the Brownian motion
$$ (-S_{i,j}B^{j}(t)+\mu t)_{\tau_i \leq t< \tau_{i+1}\wedge T};$$ that is, for  $t \in [\tau_i, \tau_{i+1}\wedge T)$ we define
\[
\tilde \Upsilon_j(t) = \tilde \Upsilon_j(\tau_i) -  S_{i,j}(B^j(t)-B^j(\tau_i))+\mu (t-\tau_i)
\]
Then we observe that before time $T$, $\tilde{\Upsilon}_j(t)$ and $\tilde{\Upsilon}_j(\tau_i)$ have the same sign, and therefore
\[
- |\tilde \Upsilon_j(t) | = -S_{i,j}\tilde{\Upsilon}_j(t) =  -|\tilde \Upsilon_j(\tau_i)| +  (B^{j}(t)-B^j(\tau_i)) - S_{i,j} \mu (t-\tau_i), \quad t \in [\tau_i, \tau_{i+1}\wedge T).
\]
Since we have $\Xi_j(\tilde{\nu},B^\mu,Q,I,\tau_i)\leq -|\tilde{\Upsilon}_j(\nu,B^\mu,Q,I,\tau_i)|$, we must also have that for $t\in [\tau_i,\tau_{i+1}\wedge T)$:
\begin{align*} 
   \Xi_j(\nu,B^\mu,Q,I,t) &=\Xi_j(\nu,B^\mu,Q,I,\tau_i)+B^{j}(t)-B^{j}(\tau_i)+\mu(t-\tau_i) \\
    &\leq -|\tilde{\Upsilon}_j(\nu,B^\mu , Q,I,\tau_i)| + B^{j}(t)-B^{j}(\tau_i) - S_{i,j}\mu(t-\tau_i) \\& =-|\tilde{\Upsilon}_j(\nu,B^\mu , Q,I,t)|
\end{align*} 
where we have used that $\mu\leq 0$ and $S_{i,j}\in \{-1,+1\}$.

It remains to describe what happens at the branching times. Recall the branching-selection operators $k,l$ from the previous proof, and let $\tilde{I}_{i+1}$ be the index such that $-|\tilde{\Upsilon}_{\tilde{I}_{i+1}}(\nu,B^\mu,Q,I,\tau_{i+1})|$ is the $I_{i+1}$\textsuperscript{th} smallest element of $-|\tilde{\Upsilon}(\nu,B^\mu,Q,I,\tau_{i+1})|$. So then if we define:
$$\tilde{\Upsilon}(\nu,B^\mu,Q,I,\tau_{i+1})=k\left(\Theta^N\left(\tilde{\Upsilon}(\nu,B^\mu,Q,I,\tau_{i+1}-)\right),\tilde{I}_{i+1}\right),$$
it immediately follows that if $\tau_{i+1}<T$, then
\begin{align*}\Xi(\tilde{\nu},B^\mu,Q,I,\tau_{i+1})=l\left(\Theta^N\left(\Xi(\tilde{\nu},B^\mu,Q,I,\tau_{i+1}-)\right),I_{i+1}\right)
\oop l\left(-\left|\Theta^N\left(\tilde{\Upsilon}(\nu,B^\mu,Q,I,\tau_{i+1}-)\right)\right|,I_{i+1}\right)\\
=k\left(\Theta^N\left(\tilde{\Upsilon}(\nu,B^\mu,Q,I,\tau_{i+1}-)\right),\tilde{I}_{i+1}\right)=\tilde{\Upsilon}(\nu,B^\mu,Q,I,\tau_{i+1}).
\end{align*}
Finally, all that remains is to observe that $\tilde{\Upsilon}(\nu,B^\mu,Q,I,t)_{t\geq 0}$ is a cadlag stochastic process with particles moving like Brownian motions with drift $\mu$. Moreover, at rate $N$, a uniformly randomly chosen particle branches, and by the definition of $k$, the particle furthest from the origin is removed. Therefore $\Theta^N(\tilde{\Upsilon}(\nu,B^\mu,Q,I,t)_{t\geq 0}$ is an $N$-Brownian bees process with drift $\mu$.
\end{proof}

\section{Association of the $N$-BBM}
Another property that we will need in the course of our proof is that $N$-BBM is associated. Association is a property of random variables which can essentially be thought of as meaning that `an FKG-type inequality holds' (see definition 1.1 of \cite{assoc}).

\begin{defn}
    Let $X$ be a random variable taking values in a partially ordered space $\C{E}$ with partial ordering $\leq_o$. A function $f:\C{E}\to \B{R}$ is called \textbf{increasing} if for $e_1,e_2\in \C{E}$, we have $e_1 \leq_o e_2 \implies f(e_1)\leq f(e_2)$. An event $A$ is called increasing if the indicator of $A$ is an increasing function. We call the random variable $X$ \textbf{associated} if we have $\B{E}[f(X)g(X)]\geq \B{E}[f(X)]\B{E}[g(X)]$ for all increasing functions $f,g$. 
\end{defn}

An immediate consequence of the definition of association is that if $h:\C{E}\to \C{E}$ is a monotonic increasing function, then $h(X)$ is also associated. This follows because for any increasing functions $f,g$, the monotonicity of $h$ implies that $f\circ h$ and $ g\circ h$ are also increasing functions, so $\B{E}[f(h(X))g(h(X))]\geq \B{E}[f(h(X))]\B{E}[g(h(X))]$.

Associated random variables were first studied by Esary, Proschan, and Walkup in \cite{assoc}, and much has been written about them in the literature since then, and we will show here that, by the construction of the $N$-BBM process, this property holds for the $N$-BBM process as well. In particular, we will prove that:

\begin{lemma} \label{nbbmAssoc} Fix $T\in (0,\infty)$. Then the $N$-BBM $(Z(t))_{t\in [0,T]}$, viewed as a $\C{D}([0,T],\B{R}^N)$-valued random variable, is associated for all $T$, under the partial order relation $\leq_o$ where $(X_t)_{t\in [0,T]} \leq_o (Y_t)_{t\in [0,T]} \iff X_t \oop Y_t \; \forall t\in [0,T]$. 
\end{lemma}

\begin{proof}
We will prove this fact by using the construction:
$$(Z^N(t))_{t\geq 0}=\Theta^N\left( \Xi\left(\nu,B, Q, I,t\right)_{t\geq 0}\right)$$  from Proposition \ref{generalCouple}.

The space $\C{D}([0,T],\B{R})$ has a natural partial ordering $\leq_T$ where for $x,y\in \C{D}([0,T],\B{R})$, we have $x\leq_T y$ if and only if $x(t)\leq y(t) \; \forall t\in [0,T]$. It is proven by Alexandre Legrand in \cite{levyAssoc} that any $1$-dimensional L\'{e}vy process when viewed as a $\C{D}([0,T],\B{R})$-valued random variable is associated under the partial ordering $\leq_T$ for any $T$. In particular, this includes the Brownian motions and Poisson process used to construct $Z$ above. The association of Brownian motion was already proven by David Barbato in \cite{barbato} under a different partial ordering, but Legrand's work extends this to the more natural partial ordering $\leq_T$. We also know that a uniform random variable on the set $\{1,2,\ldots,N\}$ is associated; Chebyshev's sum inequality states that for increasing functions $f,g:{1,\ldots,N}\to \B{R}$, $\frac1N \sum_{i=1}^N f(i)g(i) \geq \left(\frac1N \sum_{i=1}^N f(i)\right)\left(\frac1N \sum_{j=1}^N g(j)\right)$. 

Another important property of associated random variables is that association of random variables is preserved under taking products. In particular, if $E_i$ is an $\C{E}_i$-valued associated random variables with partial ordering $\leq_i$ in probability space $(\Omega_i,\C{F}_i,\B{P}_i)$ for $i\in \Gamma$, then $(E_i)_{i\in \Gamma}$ is an $\bigotimes_{i\in \Gamma}\C{E}_i$-valued associated random variable in the product probability space, with the product partial order $\leq_p$ (ie. $(a_i)_{i\in \Gamma} \leq_p (b_i)_{i\in \Gamma} \iff a_i \leq_i b_i$ for all $i\in \Gamma$). This is Lemma 2.5 in \cite{levyAssoc} for finite $\Gamma$, or Theorem 3 in \cite{barbato} for more general $\Gamma$. Therefore, as in the statement of Proposition \ref{generalCouple}, define a family of i.i.d Brownian motions $B=((B^{j}(t)){t\geq 0})_{1\leq j\leq N}$, a Poisson point process with intensity $N$, $Q=(Q(t))_{t\geq 0}$, and a family $I=(I_i)_{i\in \B{N}}$ of uniform random variables. Then $(B,Q,I)_{t\in [0,T]}=((B^{j}(t))_{1\leq j\leq N,t\in [0,T]},(Q(t))_{t\in [0,T]},(I_i)_{i\in \B{N}})$ is an associated random variable under the product partial ordering. 

Finally it remains to show that the $N$-BBM $(Z^N(t))_{t\in [0,T]}=\Theta^N(\Xi(\nu,B,Q,I,t))_{t\in [0,T]}$ is a monotonic function of the associated random variable $(B,Q,I)_{t\in [0,T]}$ and therefore is also associated. The monotonicity of $Z^N$ as a function of $B$ is obvious; increasing the Brownian motions driving the particles increases $Z^N$. It is also obvious that increasing $Q$ increases $Z^N$, as an increased $Q$ means that positive jumps of particles happen sooner. 

Finally to see that $\Xi(\nu,B,Q,I,t)_{t\in [0,T]}$ is also monotonic as a function $I$, we note that clearly by definition $l(v,i)\oop l(v,j)$ for $i\leq j$. Hence the $N$-BBM is a monotonic increasing function of the associated random variable $(B,Q,I)_{t\in [0,T]}$ and hence $(Z^N(t)_{t\in [0,T]}=\Theta^N(\Xi(\nu,B,Q,I,t))_{t\in [0,T]}$ is associated.
\end{proof}

\section{The Sub-critical Case: $|\mu|<\mu_c^N$}
Recall that for the sub-critical case of Theorem \ref{mainThem}, we want to prove the following:

\begin{theorem*}  
Let $(X^{N,\mu}(t))_{t\geq 0}$ be a one dimensional Brownian bee system with sub-critical drift $\mu<\mu_c^N$. Then $(X^{N,\mu}(nt_0))_{n\in\B{N}}$ is positive Harris recurrent for any $t_0>0$ and has a unique stationary measure $\pi^{N,\mu}$ on $\B{R}^N$ so that for any $\chi\in\B{R}^N$
$$\lim_{t\to\infty}\sup_{C}|\B{P}(X^{N,\mu}(t)\in C|X^{N,\mu}(0)=\chi)-\pi^{N,\mu}(C)|_{TV}=0.$$
\end{theorem*}

Therefore let us begin this section by defining positive Harris recurrence. We use Definition (2.2) in Athreya and Ney (\hspace{1sp}\cite{athreyaNey}), which is later used in the context of branching-selection particle systems by Berestycki et al. and Durrett \& Remenik (Proposition 6.5 in \cite{beesBBNP} and Proposition 3.1 in \cite{durrRemenik} respectively). As in \cite{beesBBNP}, this allows us to apply Theorems 4.1 and 6.1 of \cite{athreyaNey} which essentially say that the stochastic process has a unique invariant probability distribution towards which it converges.

\begin{defn} A stochastic process $(Y_n)_{n\in\B{N}}$ with state space $\C{Y}$ is called a \textbf{recurrent and strongly aperiodic Harris chain} if $\exists A\subseteq \C{Y}$ such that:
\begin{itemize}
    \item $\exists \varepsilon >0$ and probability measure $q$ on $A$ such that $\B{P}_\xi(Y_1\in C)\geq \varepsilon q(C) \; \; \forall \xi \in A, C\subseteq A$ 
    \item the hitting time of $A$ is almost surely finite for any initial $Y_0 \in \C{Y}$; that is $\B{P}_\xi(\tau_A < \infty)=1 \; \; \forall \xi \in \C{Y}$ where $\tau_A := \inf\{ n\geq 1: Y_n\in A\} $.
\end{itemize}
 Further, we call the proces \textbf{positive Harris recurrent} if $\sup_{\xi \in A}\B{E}_\xi[\tau_A]<\infty $, and \textbf{null Harris recurrent} otherwise.
\end{defn}

We begin by proving the additional condition for \textit{positive} Harris recurrence. In particular, we will first prove the follow proposition:

  \begin{propn} \label{finExpNBBM} Let $Z^N$ be an N-BBM system with drift $\mu > -\mu_c^N$. Let $\mathfrak{X}^-_{N,s}$ be the set of initial configurations of $Z^N$ in which the rightmost particle is at position $s\in \B{R}$; that is, $Z_N^N(0)=s$. Let the stopping time $T^\mu_0$ be defined by
  $$T^\mu_0:=\inf\{t\geq 0:Z^N_N(t)\geq 0\}$$
  Then there exist finite constants $a,b$ depending on $\mu$ such that $\sup_{x\in\mathfrak{X}^-_{N,-s}}\B{E}_x[T^\mu_0]\leq a+bs<\infty$.
  \end{propn}

\begin{proof}
    For the purpose of this proof, it will be helpful to consider, for $k\in \B{N}$, the $N$-BBM process $(Z^N(t))_{k\leq t\leq k+1}$ as a function of $Z^N(k)$ and $N$ independent binary branching Brownian motions (BBMs), $(\C{B}_{k}^1(t))_{0\leq t\leq 1},\ldots$,$(\C{B}_k^N(t))_{0\leq t\leq 1}$. The BBM $(\C{B}_k^i(t))_{0\leq t\leq 1}$ starts from a single particle at the origin and has constant drift $\mu$. We let $\C{N}_k^i(t)$ denote the number of particles in the BBM $\C{B}_k^i$ at time $t$, and label its particles, in the order in which they are born, by $B_k^{i,1}(t),\ldots B^{i,\C{N}_k^i(t)}_k(t)$. We also choose that $(\C{B}_k^i(t))_{0\leq t\leq 1}$ will be the BBM which drives the particle which is $i$\textsuperscript{th} smallest at time $k$, $Z_i^N(k)$

Then given $Z^N(k)$, we can construct $(Z^N(t))_{k\leq t\leq k+1}$ as a subset of
$$\bigcup_{i=1}^N \bigcup_{j=1}^{\C{N}_k^i(t)}\{Z_i^N(k)+B^{i,j}_j(t)\}=:U(t)$$
as follows. Particles in $U(t)$ are of two different types, `alive' and `ghost':
\begin{itemize}
    \item At time $k$, all particles are `alive'.
    \item When an `alive' particle branches, it branches into two `alive' particles, and simultaneously the leftmost `alive' particle (which may be one of the particles involved in the branching event) changes to `ghost' 
    \item When a `ghost' particle branches, it branches into two `ghost' particles.
\end{itemize}
Then at time $t$, the $N$-BBM $Z^N(t)$ is described by the positions of the $N$ `alive' particles in $U(t)$. For $t\in [k,k+1]$, write
\begin{align}
    Z^N(t)=\Phi\left(Z^N(k),\{(\C{B}^i_k(t))_{0\leq t\leq 1}:i=1,\ldots,N\}\right)
\end{align}
to see explicitly $Z^N(t)$ as a function of $Z^N(k)$ and $N$ independent BBMs. 

Now for $k\in \B{N}$, we will consider the following event $A_k := A^{(1)}_k \cap A^{(2)}_k \cap A^{(3)}_k \cap A^{(4)}_k$, where
\begin{align*}
A^{(1)}_k &=\{\C{N}^N_{k-1}(1)=N\}\\
A^{(2)}_k &=\{\C{N}^1_{k-1}(1)=\C{N}^2_{k-1}(1)=\ldots =\C{N}^{N-1}_{k-1}(1)=1\}\\
A^{(3)}_k &=\{\text{At each branching time }T_1,\ldots,T_{N-1}\text{ of }(\C{B}^N_{k-1}(t))_{0\leq t\leq 1}\text{ we have }B^{i,1}_{k-1}(T_j)<0\\
&\quad \quad \quad\text{for }i,j=1,\ldots,N-1\}\\
A^{(4)}_t &=\{B^{N,i}_{k-1}(T_j)>0\text{ and }B^{N,i}_{k-1}(t)<1\text{ for }i=1,\ldots,\C{N}^{N}_{k-1}(T_j),\;j=1,\ldots,N-1\text{ and }t\in[0,1]\}
\end{align*}
and define the sequence of stopping times by $\tau_0=0$ and for $i\geq 1$:
$$
\tau_i = \inf\{t\in \B{N}: t> \tau_{i-1},\; A_t\text{ occurs}\}
$$

In plain English, these are essentially the integer times when the rightmost particle of the system `regenerates' the system by becoming an ancestor of every particle alive in one unit of time. Conditions $A_k^{(3)},A_k^{(4)}$ ensure that at each branching time, the leftmost `alive' particle is one of the particles $Z^N_i(k-1)+B^{i,1}_{k-1}(T_j)$ for $i=1,\ldots, N-1$ and therefore the $N$ `alive' particles at time $k$ are the $N$ particles of the BBM $\C{B}_{k-1}^{N}$. Therefore, the particles alive at time $\tau_i$ are located at
\begin{equation}\label{pclPosns}\{Z^N_N(\tau_i-1)+B^{N,j}_{\tau_i-1}(1), j=1,\ldots,\C{N}^N_{\tau_i-1}(1)\}.\end{equation}

This construction of the $N$-BBM process then makes it explicitly clear that the event $A_t$ is independent of $Z^N(t-1)$; that is, we can verify whether $A_t$ occurs given only the driving BBMs. Therefore it follows that $\tau_{i+1}-\tau_i$ are i.i.d., and by definition of $\tau_i$, $\B{E}[\tau_{i+1}-\tau_i]\geq 1$. The definition of $A_t$ also makes it clear that the events $A_1,A_2,A_3,\ldots$ are independent and $\B{P}(A_1)=\B{P}(A_2)=\ldots$, therefore $\inf\{k\in \B{N}:A_k\text{ occurs}\}$ is a geometric random variable which bounds $\tau_1$ from above, therefore $\B{E}[\tau_{i+1}-\tau_i]<\infty$. Furthermore, since the processes $(\C{B}_{\tau_i-1}^{N}(t))_{0\leq t\leq 1}$ are i.i.d, we clearly have that $Z^N(\tau_{i+1})-Z^N(\tau_i)$ are i.i.d for $i\geq 1$ (although $Z^N(\tau_1)-Z^N(\tau_0)$ will depend on the initial condition). 

Then since $\tau_i\xrightarrow[i\to\infty]{a.s.}\infty$, we have $\lim_{i\to\infty}Z(\tau_i)/\tau_i= \lim_{t\to\infty}Z(t)/t = \mu+\mu_c^N$. So by \cite{janson} (equation (1.5)), we have that there exist finite constants $a',b'\geq 0$ such that
$$\B{E}\left[\inf\{i\in \B{N}:Z^N_N(\tau_{i+1})-Z^N_N(\tau_1)\geq R\}\right]\leq a'+b'R.$$

Now observe that $Z^N_N(\tau_1)-Z^N_N(0)$ is stochastically dominated by the maximum of a BBM started from $N$ particles. Also, since $\tau_1$ is bounded above by a geometric random variable, $\B{E}e^{\tau_1} < \infty$, therefore by the many-to-one lemma (see for example Lemma 2.1 in \cite{kim}), we have that 
$$\B{E}[Z_N^N(\tau_1)-Z^N_N(0)]<N\B{E}[e^{\tau_1}]\B{E}\left[\sup_{t<\tau_1}B_t\right]<\infty.$$ 

Thus it follows that for any $\xi \in \mathfrak{X}^-_{N,-s}$ we have
\begin{align*}\B{E}_\xi\left[\inf\{i\in \B{N}:Z_N^N(\tau_{i+1})\geq 0\}\right]&=\B{E}_\xi\left[\inf\{i\in \B{N}:Z^N_N(\tau_{i+1})-Z^N_N(\tau_1)\geq s-Z^N_N(\tau_1)+Z^N_N(0)\}\right]\\
&\leq a' + b'\left(s-\B{E}_\xi[Z^N_N(\tau_1)-Z^N_N(0)]\right).\end{align*}

The inequality follows by conditioning on $Z^N_N(\tau_1)$, and noting that $Z^N_N(\tau_{i+1})-Z^N_N(\tau_i)$ is dependent on the distribution of particles $Z^N_N(\tau_i)$ as seen from the rightmost particle $Z^N_N(\tau_i)$, but independent of $Z^N_N(\tau_i)$ itself; therefore independently of $Z^N_N(\tau_1)$, $Z^N_N(\tau_{i+1})-Z^N_N(\tau_i)$ are i.i.d. for $i\geq 1$. Thus
$$\B{E}_\xi[T_0^\mu] \leq \B{E}[\tau_1-\tau_0]\B{E}_\xi\left[\inf\{i\in \B{N}:Z_N^N(\tau_{i+1})\geq 0\}\right] \leq \B{E}[\tau_1-\tau_0]\left(a' + b'(s-\B{E}_\xi[Z^N_N(\tau_1)-Z^N_N(0)])\right)=: a+bs,$$
thus completing the proof.
\end{proof}

Note that the restriction to the set $x\in \mathfrak{X}^-_{N,-s}$ is necessary and the result is not true for general starting configurations, since we could start from a position in which the rightmost particle is arbitrarily far from the origin and therefore has arbitrarily large hitting time. We can start with an initial condition in which all particles except for the rightmost are arbitrarily far from the origin, but the proposition is essentially saying that, since the rightmost particle will eventually have $N$ descendants, the arbitrarily far left positions will not ultimately matter. Note also that the drift $-\mu_c^N$ is the critical drift at which the expected time becomes infinite. We are now ready to prove the positive Harris recurrence of the Brownian bees:

\begin{theorem} \label{BeesPosRec} Let $X^{N,\mu}$ be an $N$ Brownian bee system with drift $\mu \in(-\mu_c^N,\mu_c^N)$. Then $(X^{N,\mu}(t_0 n))_{n\in \B{N}}$ is positive Harris recurrent for any $t_0 > 0$
\end{theorem}

\begin{proof} Assume without loss of generality $\mu_c^N < \mu < 0$. In the case of $\mu=0$, Theorem \ref{mainThem} is proven in \cite{beesBBNP}. Let $A=[-1,1]^N$ and $\tau_A:=\inf\{n\geq 1:X^{N,\mu}(t_0n)\in A\} $ be the hitting time to $A$ by $(X^{N,\mu}(t_0n))_{n\in\B{N}}$. 

Define $\tilde{X}(t)$ to be the position of the particle of $X^{N,\mu}(t)$ closest to $0$. Let $\tilde{\sigma}_0:=0$ and recursively define $\tilde{\sigma}_{i+1}:=\inf\{n \in \B{N}: n > \tilde{\sigma}_{i}: \tilde{X}(s)=0\text{ for some }s\in [(n-1)t_0,nt_0]\}$. Define $\mathfrak{X}_{N,M}$ to be the set of configurations of $N$ particles with the particle closest to zero in $[-M,M]$.

Now our first observation is that as $\tilde{X}(s)=0$ for some $s \in ((\tilde{\sigma}_i-1)t_0,\tilde{\sigma}_i t_0)$, then $X^{N,\mu}((\tilde{\sigma}_i+1) t_0) \in A$ with non-zero probability bounded below by $p_0$. This can happen by the particle closest to zero branching $N-1$ times in $(\tilde{\sigma}_i t_0,(\tilde{\sigma}_i+1)t_0)$ and the Brownian motions driving the particles having sufficiently small displacements.

Next we want to control the `bad' events where $\tilde{X}(s)$ hits $0$ in $s \in ((\tilde{\sigma}_i-1)t_0,\tilde{\sigma}_i t_0)$ but we fail to have $X^{N,\mu}((\tilde{\sigma}_i+1)t_0)\in A$. Using crude bounds on Brownian motion $B$, we can see that
\begin{align*}\B{P}\Big(|\tilde{X}((\tilde{\sigma}_1+1)t_0)| &\in (M,M+1]\Big)\leq \B{P}\left(\sup_{t\leq 2t_0}|B_{2t_0}-B_t| > M\right) \\
=&\B{P}\left(\sup_{t\leq 2t_0}B_t + \mu_t > M\right)+\B{P}\left(\inf_{t\leq 2t_0}B_t + \mu_t < -M\right)=2\B{P}(\tilde{B}_{2t_0}>M-\mu t_0)+2\B{P}(\tilde{B}_{2t_0}>M+\mu t_0)\end{align*}
where $\tilde{B}$ is a standard Brownian motion and $M \in \B{N}$. This follows from the fact that $|\tilde{X}(t)|$ is bounded above by $|B(t)|$, and therefore $\B{P}(|\tilde{X}|\in (M,M+1])\leq \B{P}(|B(t)|>M)$. The final equality is an immediate application of the reflection principle. Then using the Chernoff bound, we can bound this above by
$$
\B{P}\left(|\tilde{X}((\tilde{\sigma}_1+1)t_0)| \in (M,M+1]\right)\leq \bar{C}\exp\left(-\frac{(M-\mu t_0)^2}{4t_0}\right) + \bar{C}\exp\left(-\frac{(M+\mu t_0)^2}{4t_0}\right),
$$
for some constant $\bar{C}$.

Then by Lemma \ref{absCouple}, there is a coupling between drifted $N$-BBM, $Z^{N,\mu}$, and $X^{N,\mu}$ such that before $\tilde{X}(t)$ hits zero, we have $Z^{N,\mu}(0)=-|X^{N,\mu}(0)| \implies Z^{N,\mu}(t)\oop -|X^{N,\mu}(t)|$. Thus $\tilde{X}(s)$ must hit $0$ before $Z^{N,\mu}_N(s)$ hits $0$, and we can use the coupling of Lemma \ref{absCouple} and Proposition \ref{finExpNBBM} to give the upper bound for $M\in \B{N}$
$$\sup_{x\in \mathfrak{X}_{N,M}}\B{E}_x[\tilde{\sigma}_1]<a + bM \leq C_\mu M,$$
for some constant $C_\mu<\infty$ which depends on $\mu$ and $t_0$. Therefore combining our two bounds, we show that the return time of $\tilde{X}$ to $0$ has finite expectation:
$$
\B{E}[\tilde{\sigma}_{i+1}-\tilde{\sigma}_i]\leq \sum_{M=1}^\infty \sup_{x\in\mathfrak{X}_{N,-M-1}} \B{E}_x[\tilde{\sigma}_1]\B{P}(|\tilde{X}(\tilde{\sigma}_i)|\in (M,M+1]) =: C < \infty
$$

Then noting that at each time $(\tilde{\sigma}_i+1)t_0$, $\tilde{X}((\tilde{\sigma}_i+1)t_0)\in A$ with probability at least $p_0$, and therefore $\sup_{x\in A}\B{E}_x[\tau_A]\leq C\B{E}[Geo(p_0)]<\infty$, where $Geo(p_0)$ is a geometric random variable with parameter $p_0$. 

Therefore to prove positive Harris recurrence it remains only to show that $\exists \varepsilon>0$ and measure $q$ such that 
$$\B{P}_a(X^{N,\mu}(t_0)\in S)\geq \varepsilon q(S)\quad \text{for all }a\in A, S\subseteq A.$$  

Fix $a=(a_1,...,a_N)\in A$, $S \subseteq A$. Then conditioning on the event that no branching occurs in $[0,t_0]$ (which happens with probability $e^{-t_0N}$), we have:
\begin{align} \label{recurrenceLeb}
\B{P}(X^N(t_0)\in S&|X^N(0)=a) \geq e^{-t_0N} \int_S \prod_{i=1}^N(2\pi t_0)^{-1/2}e^{-\frac{(x_i-a_i)^2}{2t_0}}dx_N...dx_1 \\
&\geq e^{-t_0N} (2\pi t_0)^{-N/2}\int_S\prod_{i=1}^Ne^{-1/t_0}dx_N...dx_1 = e^{-t_0N}(2\pi t_0)^{-N/2}e^{-N/t_0}Leb(S)
\end{align}
So normalising on $A$, which has finite Lebesgue measure, the condition is met. Thus $(X^{N,\mu}(t_0 n))_{n\in\B{N}}$ is positive Harris recurrent for $\mu \in (-\mu_c^N,\mu_c^N)$
\end{proof}

\begin{proof} \textit{(Of Theorem \ref{mainThem}.1)} We have proven that for any $t_0>0$ the process $(X^{N,\mu}(nt_0))_{n\in \B{N}}$ is a positive recurrent and strongly aperiodic Harris chain. Therefore by Theorems 4.1 and 6.1 of \cite{athreyaNey}, the conclusion of the Theorem \ref{mainThem}.1 holds. 
\end{proof}

\section{The Super-critical Case: $|\mu|>\mu_c^N$}
The coupling of Lemma \ref{generalCouple} also immediately gives proof of the transience of the $N$-Brownian bees system. Recall that for Theorem 1.3, we want to prove the following.

\begin{theorem*}
    Let $X^{N, \mu}(t)$ be a $N$-Brownian bees process with drift $|\mu| > \mu_c^N$ and initial condition $X^{N,\mu}(0)=\nu$. Then 
    $$\lim_{t \rightarrow \infty} \frac{X^{N, \mu}_1(t)}{t} = \lim_{t \rightarrow \infty} \frac{X^{N, \mu}_N(t)}{t} =  \left(|\mu| - \mu_c^N\right) \text{sign}(\mu) \text{   a.s.}$$
\end{theorem*}

The proof of this result is due to Flynn and can be found as Theorem 3.0.2 of \cite{flynn}. We give the argument in full for completeness.

\begin{proof}
    (Due to Flynn, Theorem 3.0.2 of \cite{flynn}) Suppose without loss of generality that $\mu < -\mu_c^N$. Let $T$ be the random time defined by
    $$T:=\inf\left\{s\geq 0: t\geq s \implies \frac{X_N^{N,\mu}(t)}{t}\leq \frac{\mu + \mu_c^N}{2}<0\right\}$$

    This is almost surely finite since, by the coupling of Proposition \ref{generalCouple}
    $$\liminf_{t \rightarrow \infty} \frac{X^{N, \mu}_N(t)}{t}  \leq \lim_{t \rightarrow \infty } \frac{Z^{N, \mu}_N(t)}{t} = \mu + \mu_c^N < 0 .$$
    
    Fix $\epsilon>0$. Since $T$ is almost surely finite, we can choose a $K \in \mathbb{R}$ such that $\mathbb{P}(T \geq K) \leq \epsilon$. Now we construct a process $Y^{N, \mu, K}$ such that $(Y^{N, \mu, K}(t))_{0\leq t\leq K}$ is an $N$-Brownian bees process with drift $\mu$ started from $Y^{N,\mu,K}(0) = X^{N, \mu}(0)$ and $(Y^{N, \mu, K}(t))_{t\geq K}$ is an $N$-BBM process with drift $\mu$ started from $Y^{N,\mu,K}(K)$. 
    $$\lim_{t \rightarrow \infty} \frac{Y^{N, \mu, K}_1(t)}{t} = \lim_{t \rightarrow \infty} \frac{Y^{N, \mu, K}_N(t)}{t} = \mu_c^N + \mu$$
    Moreover, by construction of $K$, with probability at least $1- \epsilon$, the processes $Y^{N, \mu, K}$ and $X^{N, \mu}$ will be identical - since if all the particles lie on the left of $0$ then the $N$ Brownian bees moves identically to the $N$-BBM, hence
    $$\mathbb{P} \left( \lim_{t \rightarrow \infty} \frac{X^{N, \mu}_1(t)}{t} = \lim_{t \rightarrow \infty} \frac{Y^{N, \mu, K}_1(t)}{t}\right) \geq 1 - \epsilon$$
    and similarly for $X^{N, \mu}_N(t)$, $Y^{N, \mu, K}_N(t)$.
    So since $\epsilon$ is arbitrary we deduce that:
    $$\lim_{t \rightarrow \infty} \frac{X^{N, \mu}_1(t)}{t} = \lim_{t \rightarrow \infty} \frac{X^{N, \mu}_N(t)}{t} = \mu + \mu_c^N \quad a.s.$$
\end{proof}

\section{The Critical Case: $|\mu|=\mu_c^N$}
Recall that we are trying to prove here that when $|\mu|=\mu_c^N$, the system satisfies an \textit{invariance principle}
$$
\left(m^{-1/2}X^{N,\mu_c^N}_i(mt)\right)_{t\geq 0} \xrightarrow[m\to\infty]{d} (\frac{\mu}{|\mu|}\beta^{-1/2}\sigma |B(t)|)_{t\geq 0}
$$
in the Skorokhod topology on $\C{D}([0,\infty),\B{R})$ for any particle $X^{N,\mu}_i$, where $B(t)$ is a standard Brownian motion. We want to further prove that the radius of the cloud of particles converges to zero in the scaling limit; that is
$$m^{-1/2}\left(X^{N,\mu_c^N}_N(mt) - X^{N,\mu_c^N}_1(mt)\right) \overset{d}{=} 0.$$

Therefore, in the scaling limit, the cloud of particles behaves like a single point mass moving about as reflected Brownian motion. We start with an analogous result for the $N$-BBM system which we will then couple to the $N$-Brownian bees. For ease of notation, for the rest of this work we will denote the $N$-BBM, $Z^{N,\mu_c^N}$, with positive critical drift $\mu_c^N$ and killing on the right simply as $Z$ and the $N$-Brownian bees $X^{N,\mu_c^N}$ with positive critical drift $\mu_c^N$ simply as $X$.

Notice that we want to prove a statement of convergence in the space $\C{D}([0,\infty),\B{R})$, the space of cadlag functions $[0,\infty)\to \B{R}$. However, in it's standard formulation, Donsker's theorem, and the machinery of the paper \cite{bbb} which we use here, are for convergence in the space $\C{D}([0,1],\B{R})$. Therefore before proceeding to the following proposition, we state some results from \cite{bills} which will be needed here. Let
\begin{align*}
g_m(t)=(m-t)\is_{\{m-1\leq t\leq m\}} + \is_{\{t\leq m-1\}}
\end{align*}

So for any $m>0$, $\psi_m:(X(t))_{t\geq 0} \mapsto (X(t)g_m(t))_{0\leq t\leq m}$ is a map $\C{D}([0,\infty),\B{R})\to \C{D}([0,m],\B{R})$ such that $X(m)g_m(m)=0$. Then

\begin{lemma} \label{DinftyConv} (Lemma 3 in Section 16 of \cite{bills}) A necessary and sufficient condition for $(Y_n(t))_{t\geq 0}$ to converge to $(Y(t))_{t\geq 0}$ in $\C{D}([0,\infty),\B{R})$ is that $\psi_m((Y_n(t))_{t\geq 0})$ converges to $\psi_m((Y(t))_{t\geq 0})$ in $\C{D}([0,m],\B{R})$ for every $m>0$
\end{lemma}

\begin{propn} \label{critNBBMInvar} Fix some $j$ in $\{1,\ldots,N\}$. Then
$$\left(m^{-1/2}Z_j(mt)\right)_{t\geq 0} \xrightarrow[m\to\infty]{d} (\beta^{-1/2}\sigma B(t))_{t\geq 0}$$
in the Skorokhod topology on $\C{D}([0,\infty),\B{R})$, where $B$ is a standard Brownian motion, and $\beta,\sigma$ are constants which we will define in the course of the proof. 
\end{propn}

\begin{proof} The first part of this proof will follow almost exactly as the proof of Proposition 1.4 in \cite{bbb}, which we will give a brief summary of now. Suppose we have a stochastic process $W$ and a sequence of finite stopping times, $\tau_1,\tau_2,\tau_3,\ldots$ satisfying the following properties
\begin{enumerate}
    \item $\tau_{i+1}-\tau_i$ are i.i.d with non-zero and finite mean for $i\geq 1$
    \item $W(\tau_{i+1})-W(\tau_i)$ are i.i.d for $i\geq 1$
    \item $\sup_{\tau_i \leq t\leq \tau_{i+1}} |W(t)-W(\tau_i)|$ are i.i.d for $i\geq 1$
    \item $\B{E}\left[\sup_{\tau_i\leq t\leq \tau_{i+1}}|W(t)-W(\tau_i)|^2\right]<\infty$
\end{enumerate}
then $W$ converges to a Brownian motion under the scaling $m^{-1/2}W(mt)$ in distribution in the Skorokhod topology on $\C{D}([0,\infty),\B{R})$. This is a simple corollary of Proposition 1.4 in \cite{bbb}, who go a step further by extending this result to the give a scaling limit for the barycentre of an $N$ particle system, however for our purposes it will be sufficient to show the four properties above.

As in Proposition \ref{finExpNBBM}, we consider the $N$-BBM process at some time $t\in [k,k+1]$, for $k\in \B{N}$ as a function of $Z(k)$ and $N$ i.i.d. BBMs. Therefore for $i=1,\ldots, N$, let $$\C{B}^i_k(t)_{0\leq t\leq 1}=\{B^{i,j}_k(t):j=1,\ldots,\C{N}^i_k(t)\}$$ be as in Proposition \ref{finExpNBBM}, but here with drift $\mu=\mu_c^N$. Let the events $A_k$ and the regeneration times $\tau_i$ also be the same as in Proposition \ref{finExpNBBM}.

Now in Proposition \ref{finExpNBBM}, we showed that $\{\tau_{i+1}-\tau_i:i\geq 1\}$ are i.i.d. random variables with finite mean, and that $\{Z(\tau_{i+1})-Z(\tau_i):i\geq 1\}$ are i.i.d. random variables. It is also immediate from construction that $\left\{\sup_{\tau_i\leq t \leq \tau_{i+1}}|Z_j(t)-Z_j(\tau_i)|:i\geq 1\right\}$ is an i.i.d family of random variables. Therefore in order to apply the result of Proposition 1.4 in \cite{bbb} it only remains to prove that
\begin{align} \label{sqBound}
    \B{E}\left[\sup_{\tau_i\leq t \leq \tau_{i+1}}|Z_j(t)-Z_j(\tau_i)|^2\right]<\infty.
\end{align}

Note that $\sup_{\tau_i \leq t\leq \tau_{i+1}}|Z_j(t)-Z_j(\tau_i)|$ is bounded above by $|Z_N(\tau_i)-Z_1(\tau_i)|$ plus $\lceil \tau_{i+1}-\tau_i\rceil$ random variables distributed like $\max_{i=1}^{\C{N}(1)}\sup_{0\leq t\leq 1}|B^i(t)|$, where $\C{B}(t)=\{B^i(t), i=1,\ldots,\C{N}(t)\}$ is a branching Brownian motion. Therefore by the Cauchy-Schwarz inequality, $\sup_{\tau_i \leq t\leq \tau_{i+1}}|Z_j(t)-Z_j(\tau_i)|^2$ is bounded above by
\begin{align}\label{boundforSq}
(\lceil \tau_{i+1}-\tau_i \rceil + 1)\left(|Z_N(\tau_i) -Z_1(\tau_i)|^2 + \Omega_1^2 + \ldots \Omega_{\lceil \tau_{i+1}-\tau_i \rceil}^2\right),
\end{align}
where the $\Omega_i$'s are i.i.d random variables distributed like $\max_{i=1}^{\C{N}(1)}\sup_{0\leq t\leq 1}|B^i(t)|$. Since the maximum of $\C{N}(1)$ positive quantities is bounded above by their sum, thus by the many-to-one lemma:
\begin{align*}
    \B{E}\left[\max_{j=1}^{\C{N}(1)} \sup_{0 \leq t\leq 1}|B^i(t)|^2\right] \leq e\B{E}\left[ \sup_{0\leq t\leq 1}|B(t)|^2\right].
\end{align*}
Then by the reflection principle:
\begin{align*}
    \B{E}\left[\sup_{0\leq t\leq 1}|B(t)|^2 \right]&=\int_0^\infty \B{P}(\sup_{0\leq t\leq 1}|B(t)|^2 > x)dx = 2\int_0^\infty \B{P}(B(1)>\sqrt{x})dx = 2\int_0^\infty \erfc(\sqrt{x})dx<\infty
\end{align*}
Therefore since $\B{E}[(\tau_{i+1}-\tau_i)^2]<\infty$ (because $\tau_{i+1}-\tau_i$ is bounded above by a geometric random variable), thus (\ref{boundforSq}) is finite, and hence (\ref{sqBound}) is finite, as required. 

Therefore our process $Z$ satisfies conditions (1)-(4) of Proposition 1.4 of \cite{bbb}. Now define $\beta=\B{E}[\tau_2-\tau_1]$ and $\sigma^2=Var(Z_j(\tau_2)-Z_j(\tau_1))$. Fix any $T\in (0,\infty)$. Then by combining the Lemmas 5.2, 5.3 and Theorems 5.1, 5.5 of \cite{bbb}, we can conclude that the following invariance principle holds:
$$
\left(m^{-1/2}(Z_j(mt)- mt\B{E}[Z_j(\tau_2)-Z_j(\tau_1)])\right)_{0\leq t\leq T} \xrightarrow[m\to\infty]{d} (\beta^{-1/2}\sigma B(t))_{0\leq t\leq T}
$$
in the Skorokhod topology on $\C{D}([0,T],\B{R})$, for standard Brownian motion $B$. Whislt \cite{bbb} states the result only for the space $\C{D}([0,1],\B{R})$; all spaces $\C{D}([0,T],\B{R})$ are essentially analogous, in the sense that the choice of finite interval $[0,T]$ is arbitrary. 

We now show that $\B{E}[Z_j(\tau_2)-Z_j(\tau_1)]=0$. Suppose for contradiction that $\B{E}[Z_j(\tau_2)-Z_j(\tau_1)]=c\neq 0$. Then $\frac{\B{E}[Z_j(\tau_i)]}{\tau_i}\to \frac{c}{\B{E}[\tau_2-\tau_1]}\neq 0$. This is a contradiction, as $Z$ has asymptotic drift $\lim_{t}Z_j(t)/t=\mu_c^N-\mu_c^N=0$. Thus:
$$
\left(m^{-1/2}Z_j(mt)\right)_{0\leq t\leq T} \xrightarrow[m\to\infty]{d} (\beta^{-1/2}\sigma B(t))_{0\leq t\leq T}.
$$
Since the map 
$$\psi_T:(m^{-1/2}Z_j(mt))_{0\leq t\leq T} \mapsto (m^{-1/2}Z_j(mt)g_T(t))_{0\leq t\leq T}$$ is continuous (see the proof of Lemma 3 in Section 16 of \cite{bills}), we can apply the mapping theorem (Theorem 2.7 of \cite{bills}), to get that
$$
\psi_T\left((m^{-1/2}Z_j(mt))_{t\geq 0}\right) \xrightarrow[m\to\infty]{d} \psi_T\left((\beta^{-1/2}\sigma B(t))_{t\geq 0}\right),
$$
and therefore by Lemma \ref{DinftyConv}, since our choice of $T$ was arbitrary, we have
$$
\left(m^{-1/2}Z_j(mt)\right)_{t\geq 0} \xrightarrow[m\to\infty]{d} (\beta^{-1/2}\sigma B(t))_{t\geq 0},
$$
as required.
\end{proof}

\begin{corollary} \label{nbbmInvarCor} $$\left(m^{-1/2}Z(mt)\right)_{t\geq 0} \xrightarrow[m\to\infty]{d} (\beta^{-1/2}\sigma B(t)\underline{1})_{t\geq 0},$$
in the Skorokhod topology, where $\underline{1}=(1,\ldots,1)\in \B{R}^N.$
\end{corollary}

\begin{proof} Let $\tau_1,\tau_2,\tau_3,\ldots$ be the regeneration times from Proposition \ref{critNBBMInvar}. Note that since $\tau_i \to \infty$ as $i\to \infty$ and $Z_N(\tau_i)-Z_1(\tau_i), i\geq 0$ are i.i.d with finite mean, thus we have $\tau_i^{-1/2}|Z_N(\tau_i)-Z_1(\tau_i)|\xrightarrow{\B{P}} 0$. And therefore in the scaling limit, the radius of the particle system goes to zero. Then by Slutsky's theorem, we have
\begin{align*}\left( m^{-1/2}Z(mt)\right)_{t\geq 0} &=\left( m^{-1/2}Z_1(mt)\underline{1}\right)_{t\geq 0}+\left(0, m^{-1/2}(Z_2(mt)-Z_1(mt)),\ldots,m^{-1/2}(Z_N(mt)-Z_1(mt)\right)_{t\geq 0}\\
&\xrightarrow[m\to\infty]{d}\left(\beta^{-1/2}\sigma B(t)\underline{1}\right)_{t\geq 0}+(0,\ldots,0)=\left(\beta^{-1/2}\sigma B(t)\underline{1}\right)_{t\geq 0}
\end{align*}
\end{proof}

Given that $Z$, in the scaling limit, converges to a Brownian motion with zero drift, we may naturally expect the properties of the Brownian motions `null-recurrence' to also hold for the $N$-BBM with critical drift, $Z$. Thus in the next Proposition we will show that the expected hitting time of $0$ by $Z_1$, uniformly over all initial configurations with $Z_1(0)=1$, is infinite. This will also be helpful for a later result. 

\begin{propn} \label{infHitNBBM} Let $\mathfrak{X}^+_{N,1}$ be the set of configurations of $N$ particles with the leftmost particle at $1$. Let $\hat{\tau}:=\inf\{t\geq 0: Z_1(t) < 0\}$ be the first hitting time of $0$ by $Z_1$. Then $\B{E}[\hat{\tau}|Z(0)=\rho]=\infty$ for all $\rho \in \mathfrak{X}^+_{N,1}$. 
\end{propn}

\begin{proof} To prove this result we will consider the process $Z_1$ as a \textit{regenerative process} (see Example 1.3 of \cite{fpMomentsPRW}). Note that $Z$ is a Markov process, but $Z_1$ is not. However, if we consider $Z_1$ specifically at the regeneration times $\tau_1,\tau_2,\ldots$ defined in Proposition \ref{finExpNBBM}, then the process $(Z_1(\tau_i))_{i\geq 1}$ is a Markov process. Therefore define for $n\geq 1$
$$\xi_n :=Z_1(\tau_{n+1}) - Z_1(\tau_n)\text{ and }\eta_n:=\inf_{\tau_{n}\leq t< \tau_{n+1}}Z_1(t)-Z_1(\tau_{n}),$$
so that $\inf_{t\geq \tau_1}Z_1(t) = \inf_{n\geq 1}\{Z_1(\tau_1) + \xi_1 + \cdots \xi_{n-1} + \eta_n\}$. The strategy of the proof will be to show that the hitting time of the barrier $n^{3/8}$ by the mean-centred random walk $S_n := \xi_1 + \cdots + \xi_n$ has infinite expectation and that $\eta_n \leq n^{3/8}$ for all $n$ with positive probability. Putting these two results together we will then therefore be able to show that the hitting time of zero by the process $(Z_1(\tau_1)+S_{n-1} + \eta_n)_{n\geq 1}$ has infinite expectation.

By the properties described in Proposition \ref{finExpNBBM}, we know that $S_n:=\xi_1 + \xi_2 + \cdots + \xi_n$ is a random walk with $\B{E}\xi_i = 0$ and $\B{E}\xi_i^2 < \infty$. Therefore by Theorem 3.2 of \cite{permantlePeres} we have that
\begin{align}\label{permantlePeresIneq} I := \inf_{n\geq 1} \sqrt{n}\B{P}(S_k \geq k^{3/8} \text{ for } 1 \leq k \leq n) >0.\end{align}

Now let us turn our attention to the $\eta_i$'s. By the properties which we proved in Proposition \ref{critNBBMInvar}, we know that $\{\eta_i : i\geq 1\}$ is an i.i.d family of positive random variables with $\B{E}[\eta_i^2]<\infty$. In fact, our proof easily extends to showing that $\B{E}[\eta_i^4]<\infty$. In particular, note that $\tau_{i+1}-\tau_i$ is bounded by a geometric random variable which has bounded fourth moment, and $\B{E}[\sup_{0\leq t\leq 1} |B(t)|^4] = 2\int_0^\infty \erfc(x^{1/4})dx < \infty$. Then just as we bound $\eta_i^2$ by (\ref{boundforSq}), we can use the Cauchy-Schwarz inequality to bound $\eta_i^4$ by
$$(\lceil \tau_{i+1} - \tau_i \rceil + 1)^3\left(|Z_N(\tau_i) - Z_1(\tau_i)|^4 + \Omega_1^4 + ... + \Omega_{\lceil \tau_{i+1} - \tau_i \rceil}^4 \right) < \infty.$$
Then as $\B{E}[\eta_i^4] < \infty$, we can prove that $\B{P}(\eta_i \geq -i^{3/8} \; \forall i\geq 1)>0 $. In particular, as the $\eta_i$'s are i.i.d negative random variables, then by Markov's inequality
\begin{align*} -\log\left(\B{P}(\eta_i \geq - i^{3/8} \; \forall i \geq 1)\right) &= -\log \left(\prod_{i\geq 1}\B{P}(\eta_i \geq -i^{3/8})\right) = - \sum_{i\geq 1}\log\left(\B{P}(\eta_i \geq - i^{3/8})\right) \\
&=-\sum_{i\geq 1}\log \left(1-\B{P}(\eta_i \leq -i^{3/8})\right) = -\sum_{i\geq 1}\log\left(1-\B{P}(\eta_i^4 \geq i^{3/2})\right) \\
&\leq -\sum_{i\geq 1}\log (1-i^{-3/2}\B{E}[\eta_1^4]) \leq -\sum_{i\geq 1} \frac{-i^{-3/2}\B{E}[\eta_1^4]}{1-i^{-3/2}\B{E}[\eta_1^4]} = \sum_{i\geq 1} \frac{\B{E}[\eta_1^4]}{i^{3/2} - \B{E}[\eta_1^4]} < \infty
\end{align*}
and therefore $\B{P}(\eta_i \geq -i^{3/8} \; \forall i\geq 1) > 0$.

Finally it will remain to show that
$$\B{P}(S_j > j^{3/8} \; \forall j\leq n, \eta_j \geq -j^{3/8} \; \forall j\leq n+1) \geq \B{P}(S_j > j^{3/8} \; \forall j\leq n)\B{P}(\eta_j \geq - j^{3/8} \; \forall j\leq n+1).$$

Note that the events $A=\{S_j > j^{3/8} \; \forall j\leq n\}$ and $B=\{\eta_j \geq -j^{3/8} \; \forall j\leq n\}$ are both increasing in the sense that if $X(t)\oop X'(t)$ for all $t\geq 0$ and the event $A$ (resp. $B$) occurs for the process $X$, then certainly the event $A$ (resp. $B$) occurs for the process $X'$. We prove in Lemma \ref{nbbmAssoc} that for any $T>0$ the $N$-BBM process $(Z(t))_{0\leq t\leq T}$ is associated and therefore by the definition of association of random variables, increasing events are positively correlated and therefore
\begin{align*}
    \B{P}(S_j > j^{3/8} \; \forall j\leq n, \eta_j \geq -j^{3/8} \; \forall j\leq n+1) &\geq \B{P}(S_j > j^{3/8} \; \forall j\leq n)\B{P}(\eta_j \geq - j^{3/8} \; \forall j\leq n+1) \\
    &\geq \frac{I}{\sqrt{n}}\B{P}(\eta_j \geq -j^{3/8}\forall j\geq 1)\end{align*}
where the second inequality follows from the fact that $\B{P}(\eta_j \geq -j^{3/8} \forall j\leq n)\geq \B{P}(\eta_j \geq j^{3/8} \forall j\geq 1)$ and (\ref{permantlePeresIneq}). 

Define $\sigma:=\inf\{n>0: S_{n-1} + \eta_n \leq 0\}$ to be the hitting time of $0$ by the perturbed random walk $(S_n + \eta_n)_{n\geq 1}$. Then:
\begin{align*}
    \B{E}\sigma &= \sum_{n\geq 1}\B{P}(\sigma \geq n) = \sum_{n\geq 1}\B{P}(S_{j-1} + \eta_j >0 \; \forall j\leq n) \\
    & \geq \sum_{n\geq 1}\B{P}(S_j > j^{3/8} \; \forall j\leq n-1, \eta_j \geq -j^{3/8} \;\forall j\leq n) \\
    & \geq \sum_{n\geq 1}\frac{I}{\sqrt{n}}\B{P}(\eta_j \geq -j^{3/8}\; \forall j\geq 1) = \infty
\end{align*}

Now suppose that $Z_1(\tau_1) > 0$. Therefore since $\B{E}\sigma = \infty$, $\B{E}[\tau_{i+1}-\tau_i]\in (0,\infty)$, and $\inf_{t\geq \tau_1}Z_1(t) = \inf_{n\geq 1}Z_1(\tau_1)+S_n + \eta_n$, thus it follows that certainly $\B{E}[\inf\{t\geq \tau_1: X_1(t)<0\}]=\infty$. We can therefore conclude that, since for every $\rho\in \mathfrak{X}_{N,1}^+$ we have $Z_1(\tau_1)>0$ with positive probability, we have that $\B{E}[\hat{\tau}|Z(0)=\rho]=\infty$ for every initial condition $\rho\in \mathfrak{X}^+_{N,1}$
\end{proof}

Next we prove that an $N$-Brownian bees system with positive critical drift asymptotically spends proportion 1 of its time with all particles in $[0,\infty)$.

\begin{lemma} \label{LRExcursions} Let $X$ be $N$-Brownian bees with critical positive drift. Then
$$ t^{-1}\int_0^t \is_{X_1(u)<0}du \xrightarrow[t\to\infty]{a.s.}0$$
\end{lemma}

\begin{proof} Let $S:=\inf\{s\geq 0: X_1(s)\leq 0\}$ be the first time that $X_1(s)$ exits $(0,\infty)$. Then observe that 
$$t^{-1}\int_0^t \is_{X_1(u)<0}du \leq t^{-1}\int_S^{S+t} \is_{X_1(u)<0}du$$
so we can assume without loss of generality that $X_1(0)\leq 0$. 

Now consider the following sequence of stopping times. Fix $\tau_0=0$ and for $i\geq 1$, recursively define $\sigma_{i}:=\inf\{t\geq \tau_{i-1}:X_1(t)\geq 1\}$ and $\tau_{i}:=\inf\{t\geq \sigma_i:X_1(t)\leq 0\}$. Define $N_t:=\max\{n:\tau_n \leq t\}$ to be the index of the most recent stopping time $\tau_i$. Note that all the time that $X_1$ spends in $(-\infty,0]$ is contained in the intervals $[\tau_{i-1},\sigma_i]$ for $i\geq 1$. So then
\begin{equation}\label{leftExcBound}t^{-1}\int_0^t \is_{X_1(u)<0}du \leq t^{-1}\sum_{i=1}^{N_t+1} (\sigma_i - \tau_{i-1}) = \frac{N_t+1}{t} \times \frac{1}{N_t +1}\sum_{i=1}^{N_t+1}(\sigma_i - \tau_{i-1}).\end{equation}

Therefore the strategy of this proof will be to show that $\B{E}_\xi[\tau_i - \sigma_i]=\infty$ uniformly for all configurations of $X(\sigma_i)$ and that $\B{E}_\xi[\sigma_i - \tau_{i-1}]<\infty$ uniformly for all configurations of $X(\tau_{i-1})$. Then by renewal theory, the first term of the right hand side of (\ref{leftExcBound}) converges to zero and the second term converges to a finite limit, so that that we have the desired convergence to zero by the algebra of limits.

Recall that $\mathfrak{X}^+_{N,1}$ is the set of configurations of $N$ particles with the leftmost particle at 1. Now, before $X_1$ hits 0, the process $X$ has the same behaviour as an $N$-BBM. Therefore by Proposition \ref{infHitNBBM}, $\B{E}[\tau_i - \sigma_i |X(\sigma_i)=\xi]=\infty$ for all $\xi \in \mathfrak{X}^+_{N,1}$, which is the set of possible configurations of $X(\sigma_i)$.

We now consider the renewal process $N'_t$ with renewal times distributed like $\tau_i-\tau_{i-1}|X(\sigma_i)=\underline{1}$. Since we can couple so that $\{\tau_i - \tau_{i-1} | X(\sigma_i)=\underline{1}\} \leq \{\tau_i - \tau_{i-1} |X(\sigma_i)=\xi\}$ for any $\xi \in \mathfrak{X}^+_{N,1}$, almost surely, therefore we can couple $(N_t)_{t\geq 0}$ and $(N'_t)_{t\geq 0}$ so that $N_t \leq N'_t$ by coupling the inter-renewal times of $N_t, N'_t$. Therefore by standard renewal theory (see Theorem 2.4.7 in \cite{durrett}), $(N_t+1)/t \leq (N'_t+1)/t \overset{\text{a.s.}}{\to} 0$. 

We now prove that $\frac{1}{N_t+1}\sum_{i=1}^{N_t+1}(\sigma_i - \tau_{i-1}) \to L<\infty$, by showing that $\sup_{\xi\in \mathfrak{X}^+_{N,0}}\B{E}_\xi[\sigma_i - \tau_{i-1}]<\infty$. To do this we will consider the times at which the particle closest to the origin is in $[0,\infty)$. For this part of the proof, it will be preferable to think of the $N$-Brownian bees in terms of their intrinsic labelling; that is, the process $V^N(t)$ instead of the functional $X^N(t)=\Theta^N(V^N(t))$. We will denote by $V$ the $N$-Brownian bees $V^N$ with critical positive drift $\mu_c^N$. Fix $i$ and define the stopping times $\tau'_{i,0},\tau'_{i,1},\tau'_{i,2},\tau'_{i,3},\ldots$ inductively by $\tau_{i,0}' = \tau_i$ and for $j\geq 0$,
$$\tau'_{i,j+1}:=\inf\{s > \tau'_{i,j}: V_{\ell_0}(s)\in [0,\infty), s-\tau_i\in \B{N}\},$$
where $\ell_0=\ell_0(s)$ is the index (under the intrinsic labeling) of the particle closest to the origin at time $s$.

As in the proof of Propositions \ref{finExpNBBM},\ref{critNBBMInvar}, we will consider our process, here the $N$-Brownian bees, as being embedded in $N$ independent branching Brownian motions. Therefore for $\ell=1,\ldots, N$, let $\C{B}^\ell_j(t)_{0\leq t\leq 1}=\{B^{\ell,m}_j(t):m=1,\ldots,\C{N}^\ell_j(t)\}$ be as in Proposition \ref{finExpNBBM}. $(\C{B}_j^\ell(t))_{0\leq t\leq 1}$ will be the BBM attached, over the interval $[\tau'_{i,j},\tau'_{i,j}+1]$, to the particle which is $\ell$\textsuperscript{th} smallest at time $\tau'_{i,j}$. As before, particles will have 2 types (`alive' and `ghost') and over the interval $[\tau'_{i,j},\tau'_{i,j}+1]$ will behave according to the following dynamics:
\begin{itemize}
    \item At time $\tau'_{i,j}$, all particles have type `alive'.
    \item When an `alive' particle branches, it branches into 2 `alive' particles, and simultaneously the `alive' particle furthest from the origin (which may be one of the two particles involved in the branching event) is changed from type `alive' to type `ghost'.
    \item When a `ghost' particle branches, it branches into 2 `ghost' particles. 
\end{itemize}
Then the $N$-Brownian bees systems is described by the set of $N$ `alive' particles. Explicitly, for $t\in [\tau_{i,j}',\tau_{i,j}'+1]$, write
$$V(t)=\Psi(V(\tau_{i,j}'),\{(\C{B}^\ell_j(t))_{0\leq t\leq 1}:\ell=1,\ldots,N\})$$
to denote that $V(t)$ is a function of $V(\tau_{i,j}')$ and the $N$ specified and independent BBMs.

Recall that at time $\tau'_{i,j}$, the particle closest to zero has index $\ell_0(\tau'_{i,j})$. Then we can now consider the event $E_j=E^{(1)}_j\cap E^{(2)}_j\cap E^{(3)}_j\cap E^{(4)}_j$, where
\begin{align*}
    E^{(1)}_j &= \{\C{N}_j^{\ell_0(\tau'_{i,j})}(1/2)=\C{N}_j^{\ell_0(\tau'_{i,j})}(1)=N, \C{N}_j^{\ell}(1)=1\text{ for }\ell\neq \ell_0(\tau'_{i,j})\}\\
    E^{(2)}_j &= \{\text{At each branching time }T_1,\ldots T_{N-1}\text{ of }(\C{B}^{\ell_0(\tau'_{i,j})}_j(t))_{0\leq t\leq 1}\text{ we have }\\
    &\quad \quad \quad \text{sign}(V_\ell(\tau'_{i,j}))B^{\ell, 1}_j(T_m - \tau'_{i,j})>1\text{ for }\ell\neq \ell_0(\tau'_{i,j})\text{ and }m=1,\ldots,N-1\}\\
    E^{(3)}_j &= \{0>B^{\ell_0(\tau'_{i,j}),\ell}_j(T_m-\tau'_{i,j})>-1\text{ and }0>B_j^{\ell_0(\tau'_{i,j}),\ell}(1/2)>-1\text{ for }m=1,\ldots,N-1\\
    &\quad \quad \quad \text{ and }\ell=1,\ldots,\C{N}^{\ell_0(\tau'_{i,j})}_j(T_m - \tau'_{i,j})\}\\
    E^{(4)}_j &= \{B^{\ell_0(\tau'_{i,j}),\ell}_j(1))>B^{\ell_0(\tau'_{i,j}),\ell}_j(1/2)+2\text{ for }\ell=1,\ldots,N\}
\end{align*}

In laymans terms, this is the event that particle $V_{\ell_0(\tau_{i,j}')}(\tau_{i,j}')$ branches $N-1$ times, all in the interval $[\tau'_{i,j},\tau'_{i,j}+1/2]$ whilst no other particle branches (the event $E^{(1)}_j$). During the interval $[\tau'_{i,j},\tau'_{i,j}+1/2]$, the particle $\ell_0(\tau_{i,j}')$ and all of it's descendants stay within $(V_{\ell_0(\tau_{i,j}')}(\tau'_{i,j})-1,V_{\ell_0(\tau_{i,j}')}(\tau'_{i,j}))\subseteq (-1,\infty)$ at branching events (the event $E^{(3)}_j$), whilst all other particles get further from the origin by distance at least $1$ at branching events (the event $E^{(2)}_j$). This means that at each branching time, the `alive' particle furthest from the origin is not among the descendants of particle $V_{\ell_0(\tau_{i,j}')}$, and therefore at time $\tau'_{i,j}+1/2$, the $N$ `alive' particles in the system are exactly the $N$ descendants of $V_{\ell_0(\tau_{i,j}')}$. Subsequently, the event $E^{(4)}_j$ ensures that all of these `alive' particles move at least $2$ units to the right, so that all the `alive' particles are in $[1,\infty)$. 

As in our construction in Proposition \ref{critNBBMInvar}, the probability of this event occurring is independent on the initial configuration of the system, therefore $\B{P}(E_j)=q_0>0$ uniformly for all $j$. 

Finally, we consider the time increments $\tau'_{i,j+1}-\tau'_{i,j}$. We want to prove that $$\B{E}[\tau'_{i,j+1}-\tau'_{i,j}]<\infty.$$

To see this, note that if at time $t_0$, we have $V_{\ell_0(t_0)}(t_0)<0$, then whilst $V_{\ell_0(t_0)}(t)<0$, $V_{\ell_0(t_0)}(t)$ is bounded below by a Brownian motion with drift $\mu_c^N$. This is because particles, under the intrinsic labelling, only ever jump closer to the origin. As a Brownian motion with positive drift exits $(-\infty,0]$ in finite expected time, certainly the expected time it takes until $V_{\ell_0(t)}(t)$ exits $(-\infty,0]$ is also finite. Therefore $\B{E}_\xi[\tau'_{i,j+1}-\tau'_{i,j}]<M<\infty$ uniformly for any starting configuration $V(\tau'_{i,j})=\xi\in \mathfrak{X}^+_{N,0}$. Therefore uniformly over initial configurations $X(\tau_{i})=\xi\in \mathfrak{X}^+_{N,0}$, we have $ \B{E}_\xi[\sigma_{i+1} - \tau_{i}] \leq M\B{E}[Geom(q_0)]<\infty$. Therefore by the strong law of large numbers
$$\lim_{t\to\infty}\frac{1}{N_t+1}\sum_{i=1}^{N_t+1}(\sigma_{i+1} - \tau_i)\leq M\B{E}[Geom(q_0)] < \infty \quad a.s.,$$
Hence we can conclude that the right hand side of equation (\ref{leftExcBound}) almost surely converges to $0$ as $t\to\infty$, and thus the claim holds. 
\end{proof}

Next we state a technical result about the almost sure continuity (with respect to Brownian motion) of a transformation of path space. In particular, the transformation is one which takes a path $[0,\infty)\to \B{R}$ and transforms it into a path $[0,\infty)\to [0,\infty)$ by removing the excursions in the negative half-line. Specifically, this is the transformation $g$ defined by:
$$g:(Z(t))_{t\geq 0}\mapsto \left( Z(\inf\{s\geq 0: \int_0^s \is_{Z(u)\geq 0}du \geq t\})\right)_{t\geq 0}$$

The motivation behind this technical result is that it will allow us to apply the mapping theorem (Theorem 2.7 in \cite{bills}), which states that if $A_n \to A$ in distribution and the transformation $g$ has a discontinuity set $\C{D}_g$ such that $\B{P}(A\in \C{D}_g)=0$, then $g(A_n)\to g(A)$ in distribution. Since our random variables are cadlag paths, we must prove that our transformation is almost surely continuous under the Skorokhod metric on $\C{D}([0,\infty),\B{R})$. The transformation described above is chosen because it transforms a Brownian motion into a reflected Brownian motion; therefore by the mapping theorem and Proposition \ref{critNBBMInvar}, the continuity result will give that:
$$g((m^{-1/2}Z_i(mt))_{t\geq 0}) \xrightarrow[m\to\infty]{d} (\beta^{-1/2}\sigma|B(t)|)_{t\geq 0}$$
in the Skorokhod topology on $\C{D}([0,\infty),\B{R})$ for any $i=1,2,\ldots,N$. 

Before stating the next proposition, we recall the definition of the Skorokhod metric $d_T$ on the space $\C{D}([0,T],\B{R})$ for $T>0$ and the Skorokhod metric $d_\infty$ on $\C{D}([0,\infty),\B{R})$. 

\begin{defn} Fix $T>0$ and let $\Lambda_T$ be the set of continuous and strictly increasing bijections $[0,T]\to[0,T]$. Then for $X,Y\in \C{D}([0,T],\B{R})$
\begin{align*}
    d_T(X,Y)&:=\inf_{\lambda \in \Lambda_T}\left\{\sup_{t\in [0,T]}|t-\lambda(t)|\vee \sup_{t\in [0,T]}|X(t)-Y(\lambda(T))|\right\} 
\end{align*}
and for $X,Y\in \C{D}([0,\infty),\B{R})$
\begin{align*}
    d_\infty(X,Y)&:=\sum_{T\in\B{N}}^\infty 2^{-T}(1\wedge d_T(\psi_T(X),\psi_T(Y))).
\end{align*}
\end{defn}

We can now define precisely our transformation and state the following proposition

\begin{propn} \label{pathSpaceTfm}
    Let $(W(t))_{t\geq 0}$ be a stochastic process whose sample paths are in $\C{D}([0,\infty),\B{R})$, the space of cadlag paths $[0,\infty)\to \B{R}$. Let $g:\C{D}([0,\infty),\B{R})\to\C{D}([0,\infty),\B{R})$ be the map given by:
    $$
    g:(W(t))_{t\geq 0} \mapsto \left(W(\inf\{s\geq 0: \int_0^s \is_{W(u) \geq 0} du \geq t\})\right)_{t\geq 0}.
    $$
    Let $D_g\in \C{D}([0,\infty),\B{R})$ be the discontinuity set of $g$ under the metric $d_\infty$. Then $\B{P}((B_t)_{t\geq 0} \in D_g)=0$, where $B$ is a standard Brownian motion. 
\end{propn}

The final result we need before we can prove Theorem \ref{mainThem}.2 is the following, which essentially states that changing a stochastic process by an asymptotically small time change doesn't change the scaling limit. In particular:

\begin{propn} \label{asympSmallTC} Let $(W(t))_{t\geq 0}$ be a stochastic process such that $(m^{-1/2}W(mt))_{t\geq 0}$ converges in distribution as $m\to\infty$ to $(B(t))_{t\geq 0}$ (resp. $(|B(t)|)_{t\geq 0}$) in the Skorokhod topology on $\C{D}([0,\infty),\B{R})$, where $B$ is a standard  Brownian motion. Suppose that for every $T>0$ we have
$$\sup_{0\leq t\leq T}|m^{-1}\alpha(mt)|\xrightarrow[m\to\infty]{\B{P}} 0.$$
Then $\left(m^{-1/2}W(mt+\alpha(mt))\right)_{t\geq 0}\xrightarrow[m\to\infty]{d}(B(t))_{t\geq 0}$ (resp. $(|B(t)|)_{t\geq 0}$) in the Skorohod topology on $\C{D}([0,\infty),\B{R})$
\end{propn}

The proofs of these two results are technical, and so we postpone them until section \ref{techSec}. We are now ready to prove an invariance principle to a reflected Brownian motion, which was the second regime in Theorem \ref{mainThem}. 

\begin{theorem} Let $X$ be the Brownian bees system with critical drift. Then we have the invariance principle:
$$
\left(m^{-1/2}X_1(mt)\right)_{t\geq 0} \xrightarrow[m\to\infty]{d} \left(\beta^{-1/2}\sigma|B(t)|\right)_{t\geq 0},
$$
in the Skorokhod topology on $\C{D}([0,\infty),\B{R})$, for constants $\beta,\sigma$ defined in Proposition \ref{critNBBMInvar}.
\end{theorem}

\begin{proof}
The idea of this proof is the following: since the Brownian bees model with critical drift spends asymptotically zero time to the left of the origin, it can be coupled (up to an asymptotically small time shift) to the image by a transformation similar to $g$ of an $N$-BBM process with critical drift. The invariance principle for the $N$-BBM $Z$, and the mapping theorem with the almost surely continuous map $g$ then yield the conclusion.

Recall that $Z$ is the $N$-BBM with killing on the right and critical drift of $\mu_c^N$. Let $\tilde{Z}$ be a time-change of $Z$ in which we remove the excursions where $Z_1<0$; that is
$$(\tilde{Z}(t))_{t\geq 0}:=\left(Z(\inf\{s\geq 0: \int_0^s \is_{Z_1(u)\geq 0}du \geq t\})\right)_{t\geq 0}.$$

Let $\tilde{Z}_1$ denote the position of the leftmost particle of $\tilde{Z}$ at time $t$ and note that $\tilde{Z}_1=g(Z_1)$. Therefore by the mapping theorem (Theorem 2.7 in \cite{bills}) and Proposition \ref{pathSpaceTfm}
\begin{align}\label{ztildeConv}\lim_{m\to\infty}(m^{-1/2}\tilde{Z}_1(mt))_{t\geq 0} =g\left(\lim_{m\to\infty}(m^{-1/2}Z_1(mt))_{t\geq 0}\right)=g((\beta^{-1/2}\sigma B(t))_{t\geq 0}) \overset{d}{=} (|\beta^{-1/2}\sigma B(t)|)_{t\geq 0}.\end{align}

We will now carefully construct regeneration times of $X$, by defining events $C_k$ for $k\in \B{N}$ in which the leftmost particle of $X$ is positive and, in some suitably controlled way, `regenerates' the system by becoming the parent of every particle. As with the proofs of Propositions \ref{finExpNBBM}, \ref{critNBBMInvar} and Lemma \ref{LRExcursions}, we will describe the regeneration event by considering the $N$-Brownian bees system as being embedded in $N$ independent $N$-BBMs. 

Therefore for $i=1,\ldots,N$, $k=1,2,\ldots$, and $t\geq 0$, let $\C{B}^i_k(t)=\{B^{i,j}_k(t),j=1,\ldots,\C{N}_k^i(t)\}$  be as in Proposition \ref{LRExcursions}, with the same dynamics of `alive' and `ghost' particles describing the behaviour of $X$; that is, defined so that at any time $t\in [k,k+1]$, $X(t)$ is given by the positions of the $N$ `alive' particles in the system and can be explicitly writen as
$$X(t)=\Psi(X(k),\{(\C{B}^i_k(t))_{0\leq t\leq 1}:i=1,\ldots,N\}).$$ 

Again, let the branching Brownian motion $(\C{B}_k^i(t))_{0\leq t\leq 1}$ drive the particle which is the $i$\textsuperscript{th} largest at time $k$. Now we can describe the `regeneration' events, $C_k$, of the process $X$, in terms of the branching Brownian motions and the position of the leftmost particle. Let $C_k:=C_k^{(1)}\cap C_k^{(2)}\cap C_k^{(3)}\cap C_k^{(4)}\cap C_k^{(5)}\cap C_k^{(6)}$, where
\begin{align*}
    C^{(1)}_k &= \{X_1(k-1)\geq 1\}, \\
    C^{(2)}_k &= \{\C{N}_{k-1}^1(1)=N\}, \\
    C^{(3)}_k &= \{\C{N}_{k-1}^j(1)=1\text{ for }j\neq 1\}, \\
    C^{(4)}_k &= \{\text{At each branching time $T_1,\ldots,T_{N-1}$ of $(\C{B}^1_{k-1}(t))_{0\leq t\leq 1}$, we have $B^{1,j}_{k-1}(T_\ell)<0$ for }\\
    &\quad \quad \quad \text{for $\ell=1,\ldots,N-1$ and $j=1,\ldots, \C{N}^1_{k-1}(T_\ell)=\ell+1$}\}, \\
    C^{(5)}_k &= \{B^{i,1}_{k-1}(T_\ell)>0\text{ for }\ell=1,\ldots,N-1\text{ and }i\neq 1\}, \\
    C^{(6)}_k &= \{\min_{1\leq i\leq N}\min_{1\leq j\leq \C{N}^i_{k-1}(1)}\inf_{t\in [0,1]}B^{i,j}_{k-1}(t)>-1/2 \}, 
\end{align*}
and define the regeneration times of $X$ as $\tau_0^X:=0$ and subsequently
$$\tau_{i+1}^X:=\inf\{k\in \B{N}:k>\tau_i^X,\text{ and }C_k\text{ occurs}\}.$$

In layman's terms, the event $C_k$ is an event in which the leftmost particle of the system, whilst to the right of the origin, becomes ancestor to every particle in the system over $1$ time unit. Specifically, we start with the leftmost particle in $[1,\infty)$ at time $k-1$ (the event $C^{(1)}_k$) and ask that it branches $N-1$ times in $[k-1,k]$ (the event $C^{(2)}_k$) whilst no other particle branches (the event $C^{(3)}_k$). Then the events $C^{(4)}_k$ and $C^{(6)}_k$ ensure that at each branching time $T_\ell$ of the BBM $\C{B}^1_k$, all descendants of the particle $X_1(k-1)$ are in the interval $[X_1(k-1)-1,X_1(k-1)]\subseteq [1/2,\infty)$, and that at all times \textit{all} particles remain in $[1/2,\infty)$. The event $C_k^{(5)}$ ensures that particles not descended from $X_1(k-1)$, at the branching times, are to the right of their initial position. Together, these conditions ensure that at each branching time, the particle furthest from the origin is not a descendant of the leftmost particle $X_1(k-1)$, and so at time $k$, conditional on the event $C_k$, the $N$ particles alive in the system are all descendants of particle $X_1(k-1)$.  

We now make the following observation, which demonstrates why we defined the event $C_{k}$ as we did. The events $C^{(2)}_k,\ldots,C^{(6)}_k$ are all independent of the initial condition $X(k-1)$ and depend only on the BBMs driving the system. Therefore the function $\is_{C_k}$ is a function only of $X_1(k-1)$ and the BBMs $\C{B}^1_{k-1},\ldots,\C{B}^N_{k-1}$. 

Note that, as we see in the proof of Proposition \ref{LRExcursions}, the expected time between successive events such that $\{X_1(t)\geq 1\}$ is finite, say with mean bounded above by $M$ uniformly for all configurations. Furthermore, at each time $k\in \B{N}$ such that $\{X_1(k-1)\geq 1\}$, the event $C_k$ occurs with a non-zero probability $p_0>0$, uniformly for all configurations $X(k-1)$. The number of occasions that $X_1(k-1)\geq 1$ until $C_k$ occurs is therefore a geometric random variable with mean $p_0^{-1}$, therefore certainly we have $\B{E}_\xi[\tau_{i+1}^X-\tau_i^X]\leq M\B{E}[Geom(p_0)]$ for all possible configurations $X(\tau_i^X)=\xi$. Given the process $(X(t))_{t\geq 0}$ we will now construct the process $Z$ (and by extension $\tilde{Z}$) from $X$ so that $Z$ is indeed an $N$-BBM and is coupled to $X$.

\textbf{Case 1:} If $X_1(t)\geq 0$ for all $t\in [\tau^X_0,\tau^X_1]$, then $X$ behaves exactly as an $N$-BBM process staying always above zero. Therefore define $(Z(t))_{\tau_0^X \leq t\leq \tau_1^X}:=(X(t))_{\tau_0^X \leq t\leq \tau_1^X}$, and $\tau^Z_1:=\tau_1^X$. Therefore certainly $X(\tau_1^X)=\tilde{Z}(\tau_1^X)$ and $Z$ behaves like an $N$-BBM on $[\tau_0^X,\tau_1^X]$

\textbf{Case 2:} If $X_1$ hits $0$, say for the first time at $\theta_1 \in [\tau^X_0,\tau^X_1]$, then define $(Z(t))_{\tau_0^X \leq t\leq \theta_1}:=(X(t))_{\tau_0^X \leq t\leq \theta_1}$. As all particles are above the origin in this interval, $Z$ certainly behaves like an $N$-BBM. Then on $[\theta_1,\tau_1^X-1]$, we can couple $X$ to an $N$-BBM process $Z$ as in Proposition \ref{generalCouple}, so that $Z(\theta_1)=X(\theta_1)$ and $Z(t)\oop X(t)$ for $t\in [\theta_1,\tau^X_1-1]$. Note that by definition of $\tau_1^X$, $X_1$ must stay positive in $[\tau_1^X-1,\tau_1^X]$, so certainly $\theta_1 < \tau_1^X -1$.

Now we consider the position $X_1(\tau_1^X-1)$, the position at which $X_1$ begins the next regeneration event. By our coupling (Proposition \ref{generalCouple}), we certainly have $Z_1(\tau^X_1-1)\oop  X_1(\tau^X_1-1)$. So now we continue to run the process $Z$ as an independent $N$-BBM until the time $\theta_2$ at which $Z_1(\theta_2)=X_1(\tau_1^X-1)$. Finally, for $t\in [\theta_2,\theta_2+1]$, we define $Z(t)=\Psi(Z(\theta_2),\{(B^i_{\theta_2}(t))_{0\leq t\leq 1}:i=1,\ldots,N\})$, and
$$\tau_1^Z=\inf\left\{s\geq 0: \int_0^s\is_{Z_1(u)\geq 0}du \geq \theta_2+1\right\}.$$

By the construction of our event $C_{\tau_1^X}$, all particles of $X$ stay above $0$ in $[\tau_1^X-1,\tau_1^X]$, and therefore all particles of $Z$ stay above $0$ in $[\theta_2,\theta_2+1]$. Therefore $(Z(t))_{t\in (\theta_2,\theta_2+1]}$ is exactly $(\tilde{Z}(t))_{t\in (\tau_1^Z-1,\tau_1^Z]}$, and $Z$ behaves exactly as an $N$-BBM process.

First we claim that $\B{E}_\xi|\tau_1^X - \tau_1^Z|$ is bounded uniformly for all initial conditions $X(\tau_0^X)=\tilde{Z}(\tau_0^Z)=\xi$. As above, we know that $\B{E}_\xi[\tau_{i+1}^X - \tau_i^X]\leq M\B{E}[Geom(p_0)]$. Then the boundedness of $\B{E}_\xi[\tau_1^Z]$ follows because, in case 2, $\tau_1^Z$ is at most $\theta_1$ plus the hitting time of $X_1(\tau_1^X -1)>0$ by $\tilde{Z}_1$, plus $1$. $[0,X_1(\tau_1^X-1)]$ is a compact interval, and therefore the expected hitting time of $X_1(\tau_1^X-1)$ by $\tilde{Z}_1$ is finite. Therefore $\B{E}_\xi|\tau_1^X - \tau_1^Z|$ is uniformly bounded for all $\xi$.

Next we claim that $X(\tau_1^X)=\tilde{Z}(\tau_1^Z)$. Since the event $C_{\tau_1^X-1}$ is dependent only on the driving BBMs in $[\tau_1^X-1,\tau_1^X]$ and the position $X_1(\tau_1^X-1)=\tilde{Z}(\tau^Z_1-1)$, therefore our coupling ensures that:
$$X(\tau^X_1)=\Psi(X(\tau_1^X-1),\{(\C{B}^i_{\tau_1^X-1}(t))_{0\leq t\leq 1}:i=1,\ldots,N\}) = \Psi(\tilde{Z}(\tau_1^Z-1),\{(\C{B}^i_{\tau_1^Z-1}(t))_{0\leq t\leq 1}:i=1,\ldots,N\}) = \tilde{Z}(\tau^Z_1).$$

Continuing in the same way inductively, we can define stopping times $\tau_i^X$ and $\tau_i^Z$ for $i\geq 2$ so that $i\geq 0$
\begin{align}\label{xzTildeCoupling}
X(\tau_i^X)=\tilde{Z}(\tau_i^Z) = \tilde{Z}\left(\tau_i^X + \sum_{j=1}^i ((\tau_j^Z-\tau_{j-1}^Z)-(\tau_j^X-\tau_{j-1}^X))\right).
\end{align}
and that, whenever $X_1$ exits $(0,\infty)$ in $[\tau_{j-1}^X,\tau_j^X]$, we have $\B{E}_\xi\left[|(\tau_j^Z-\tau_{j-1}^Z)-(\tau_j^X-\tau_{j-1}^X)|\right]<E_{\text{bound}}<\infty$ uniformly for all initial conditions $X(\tau_{j-1}^X)=\tilde{Z}(\tau_{j-1}^Z)=\xi$.

Now let $\tau^X(mt)$ be the smallest regeneration time $\tau_i^X$ after $mt$, and let $I(mt)$ denote its index. Let 
$$\alpha(mt):=\sum_{i=1}^{I(mt)}\left((\tau_i^Z-\tau_{i-1}^Z) - (\tau_i^X-\tau_{i-1}^X)\right),$$
so that we can write $X(\tau^X(mt))=\tilde{Z}(\tau^X(mt)+\alpha(mt))$. Fix $T\in (0,\infty)$. Using the same ideas as in the proof of Lemma \ref{LRExcursions}, we will now show that $\sup_{0\leq t\leq T}|m^{-1}\alpha(mt)|\xrightarrow[m\to\infty]{\B{P}}0.$
Define $N_{mt}:=\sum_{j=1}^{I(mt)}\is_{(\tau_j^Z-\tau_{j-1}^Z)\neq (\tau_j^X-\tau_{j-1}^X)}$. So $N_{mt}$ counts the number of intervals $[\tau_{i-1}^X,\tau_i^X]$ up to $[\tau^X_{I(mt)-1},\tau^X_{I(mt)}]$ in which $X_1$ exits $(0,\infty)$, and hence in which $(\tau^Z_j-\tau^Z_{j-1})$ and $(\tau^X_j-\tau^X_{j-1})$ are not necessarily equal. Note that at each time $\tau_i^X$, we have $X_1(\tau_i^X)>1/2$, so we can bound $N_{mt}$ above by a renewal process whose inter-arrival times are distributed like the hitting time of $0$ started from the initial configuration $(1/2,\ldots,1/2)$. By Proposition \ref{infHitNBBM}, these inter-arrival times have infinite expectation, so by standard renewal theory, $N_{mt}/mt\xrightarrow[t\to\infty]{a.s.} 0$.

Combining this with the fact above that $\B{E}_\xi[|(\tau_j^Z - \tau_{j-1}^Z) - (\tau_j^X-\tau_{j-1}^X)|]<E_{\text{bound}}$, the strong law of large numbers gives that
$$\lim_{t\to\infty}\frac{1}{N_{mt}}\sum_{i=1}^{I(mt)}\left|(\tau_i^Z-\tau_{i-1}^Z) - (\tau_i^X-\tau_{i-1}^X)\right|<E_{\text{bound}}\quad a.s.,$$
Therefore putting these convergence results together
\begin{align*}\sup_{0\leq t\leq T}|m^{-1}\alpha(mt)|\leq  \frac{N_{mt}}{m} \times \frac{1}{N_{mt}}\sum_{i=1}^{I(mt)}\left|(\tau_i^Z-\tau_{i-1}^Z) - (\tau_i^X-\tau_{i-1}^X)\right| \xrightarrow[t\to\infty]{a.s.} 0,\end{align*}
and thus, since our choice of $T$ was arbitrary, therefore by Proposition \ref{asympSmallTC}, we have
\begin{align} \label{removingAlpha}
\lim_{m\to\infty}(m^{-1/2}\tilde{Z}_1(\tau^X(mt)+\alpha(mt)))_{t\geq 0}=\lim_{m\to\infty}(m^{-1/2}\tilde{Z}_1(\tau^X(mt))_{t\geq 0}.\end{align}

Finally, we note that $\tau^X(mt)-mt$ is, in the language or renewal theory, an `excess lifetime process' or `residual waiting time' process, so it follows that for any $T\in (0,\infty)$, we have 
$$\sup_{0\leq t\leq T}\left|m^{-1}(\tau^X(mt)-mt)\right|\xrightarrow[m\to\infty]{\B{P}}0$$ (see Example 4.4.8 in \cite{durrett}). Therefore applying Proposition \ref{asympSmallTC} again, we have
\begin{align*}
\lim_{m\to\infty}(m^{-1/2}X_1(mt))_{t\geq 0}=&\lim_{m\to\infty}(m^{-1/2}X_1(\tau^X(mt)))_{t\geq 0}\quad \quad && \text{(by Proposition \ref{asympSmallTC})}\\
=&\lim_{m\to\infty}(m^{-1/2}\tilde{Z}_1(\tau^X(mt)+\alpha(mt)))_{t\geq 0}\quad \quad && \text{(by (\ref{xzTildeCoupling}))} \\
=&\lim_{m\to\infty}(m^{-1/2}\tilde{Z}_1(\tau^X(mt)))_{t\geq 0}\quad \quad && \text{(by (\ref{removingAlpha}))}\\
=&\lim_{m\to\infty}(m^{-1/2}\tilde{Z}_1(mt))\quad \quad && \text{(by Proposition \ref{asympSmallTC})}\\
\overset{d}{=}&\left(\beta^{-1/2}\sigma |B(t)|\right)_{t\geq 0}\quad \quad && \text{(by (\ref{ztildeConv}))}
\end{align*}

\end{proof}

We can now give our proof of Theorem \ref{mainThem}.2, that is, that $(X(t))_{t\geq 0}$ satisfies the invariance principle that
$$\left(m^{-1/2}X(mt)\right)_{t\geq 0}\xrightarrow[m\to\infty]{d}\left(\beta^{-1/2}\sigma|B(t)|\underline{1}\right)_{t\geq 0},$$ where $B$ is a standard Brownian motion and $\underline{1}=(1,\ldots,1)\in \B{R}^N$.

\begin{proof} (of Theorem 1.2) Let $\tau_1^X,\tau_2^X,\tau_3^X,\ldots$ be defined as in the above theorem. Note that since $\tau_i^X > i$, thus $\tau_i^X \to \infty$ as $i\to \infty$. Furthermore, the random variables $X_N(\tau^X_i)-X_1(\tau^X_i), i\geq 0$ are i.i.d with finite mean, thus we have $\tau_i^{-1/2}|X_N(\tau^X_i)-X_1(\tau^X_i)|\xrightarrow[i\to\infty]{\B{P}} 0$. And therefore in the scaling limit, the radius of the particle system goes to zero. Then by Slutsky's theorem, we have
\begin{align*}\left( m^{-1/2}X(mt)\right)_{t\geq 0} &=\left( m^{-1/2}X_1(mt)\underline{1}\right)_{t\geq 0}+\left(0, m^{-1/2}(X_2(mt)-X_1(mt)),\ldots,m^{-1/2}(X_N(mt)-X_1(mt))\right)_{t\geq 0}\\
&\xrightarrow[m\to\infty]{d}\left(\beta^{-1/2}\sigma B(t)\underline{1}\right)_{t\geq 0}+(0,\ldots,0)=\left(\beta^{-1/2}\sigma B(t)\underline{1}\right)_{t\geq 0}
\end{align*}
\end{proof}


\section{Proof of Proposition \ref{pathSpaceTfm} and Proposition \ref{asympSmallTC}} \label{techSec}
Recall Proposition \ref{pathSpaceTfm}:
\begin{propn-non}
    Let $(W(t))_{t\geq 0}$ be a stochastic process whose sample paths are in $\C{D}([0,\infty),\B{R})$, the space of cadlag paths $[0,\infty)\to \B{R}$. Let $g:\C{D}([0,\infty),\B{R})\to\C{D}([0,\infty),\B{R})$ be the map given by:
    $$
    g:(W(t))_{t\geq 0} \mapsto \left(W(\inf\{s\geq 0: \int_0^s \is_{W(u) \geq 0} du \geq t\})\right)_{t\geq 0}.
    $$
    Let $D_g\in \C{D}([0,\infty),\B{R})$ be the discontinuity set of $g$ under the metric $d_\infty$. Then $\B{P}((B_t)_{t\geq 0} \in D_g)=0$, where $B$ is a standard Brownian motion. 
\end{propn-non}

\begin{proof}
    Let $Y$ be a path in $\C{D}([0,\infty),\B{R})$ such that $\{t:Y(t)=0\}$ is closed and Lebesgue null, and $Y$ is $\alpha$-H\"{o}lder continuous for $\alpha \in (0,1/2)$. Note that Brownian motion satisfies these properties almost surely. Fix $T\in (0,\infty)$ and $\epsilon>0$. We will first show that given $\delta>0$ sufficiently small, $d_\infty(X,Y)<\delta$ implies that $d_T(g(X),g(Y))<\epsilon$. In a slight abuse of notation, we use $d_T(g(X),g(Y))$ to mean the metric $d_T$ on $\C{D}([0,T],\B{R})$ of the paths $g(X)$ and $g(Y)$ \textit{restricted to the interval} $[0,T]$.

    Let $\tilde{T}$ be the smallest integer time such that $\int_0^{\tilde{T}}\is_{Y(u)\geq 0}du \geq T+1$. Since $d_\infty(X,Y)<\delta$, we have $d_{\tilde{T}}(X,Y)<2^T\delta=:\tilde{\delta}$. So there exists a continuous and strictly increasing bijection $\lambda:[0,\tilde{T}]\to [0,\tilde{T}]$ such that $\sup_{t\in [0,\tilde{T}]}|t-\lambda(t)|<\tilde{\delta}$ and $\sup_{t\in [0,\tilde{T}]}|X(t)-Y(\lambda(t))|<\tilde{\delta}$. We will now show that by choosing $\tilde{\delta}$ sufficiently small, we can make $|\int_0^s \is_{X(u)\geq 0}du-\int_0^s \is_{Y(\lambda(u))\geq 0}du|$ arbitrarily small for all $s\in [0,\tilde{T}]$. 

    Now note that the zero set of $Y$ in $[0,\tilde{T}]$, $S_0:=\{t\in [0,\tilde{T}]:Y(t)=0\}\subseteq [0,\tilde{T}]$ is closed and bounded (and therefore compact) and Lebesgue null. Therefore the zero set of $Y\circ \lambda$, $S^\lambda_0:=\{t\in[0,\tilde{T}]:Y(\lambda(t))=0\}$, which is the image of $S_0$ under the continuous bijection $\lambda$, is also closed. Now let $A:=[0,\tilde{T}]\setminus S^\lambda_0$ be the complement of $S^\lambda_0$, which is thus open, and hence can be written as a countable union of disjoint intervals, $A=\bigcup_{i=1}^\infty A_i$. 
    
     Now since $\lambda^{-1}$ is a bijection, $\lambda^{-1}(A_i)$ are intervals contained in $[0,\tilde{T}]\setminus S_0$, the complement of $S_0$. And since $\sup_t|t-\lambda(t)| < \tilde{\delta}$, we have that $Leb(\lambda^{-1}(A_i))\in (Leb(A_i)-\tilde{\delta},Leb(A_i)+\tilde{\delta})$. Therefore choosing $k$ sufficiently large, $Leb\left(\cup_{i=1}^k \lambda^{-1}(A_i)\right)>\tilde{T}-\epsilon/4$. Then choosing $\tilde{\delta} < \epsilon/8k$, we have that $Leb\left(\cup_{i=1}^k A_i\right)>Leb\left(\cup_{i=1}^k \lambda^{-1}(A_i)\right)-2k\tilde{\delta} > \tilde{T}-\epsilon/2$. Then given an open interval $A_i=(a_i,b_i)$, define the closed interval $I_i=[a_i+\frac{\epsilon}{4k},b_i-\frac{\epsilon}{4k}]$, so that $\bigcup_{i=1}^k I_i$ is a finite set of closed intervals contained in $A$ such that $Leb\left( \bigcup_{i=1}^k I_i\right)>Leb(\bigcup_{i=1}^k A_i)-2k\frac{\epsilon}{4k} > \tilde{T}-\epsilon$.

    By the compactness of $\bigcup_{i=1}^k I_i$, the infimum $\inf\{|Y(\lambda(t))|:t\in \cup_{i=1}^k I_i\}$ is attained in $\bigcup_{i=1}^k I_i$ and thus is $>0$. Therefore choosing $\tilde{\delta} < \inf\{|Y(\lambda(t))|:t\in \cup_{i=1}^k I_i\}$, the condition $\sup_{t\leq \tilde{T}}|X(t)-Y(\lambda(t))|<\tilde{\delta}$ ensures that on $\bigcup I_i$ we have either $X(t),Y(\lambda(t))> 0$ or $X(t),Y(\lambda(t))<0$. Therefore for any $s\in [0,\tilde{T}]$
    \begin{align}\label{firstLambdaBound} \left|\int_0^s \is_{X(u)\geq 0}du - \int_0^s \is_{Y(\lambda(u))\geq 0}du\right| \leq \int_0^s |\is_{X(t)\geq 0}-\is_{Y(\lambda(t))\geq 0}|du \leq \int_0^{\tilde{T}} \is_{\{u\notin \bigcup I_i\}}du = \tilde{T}-Leb \left(\bigcup I_i\right) < \epsilon.\end{align}

    Similarly, since $Y$ is $\alpha$-H\"{o}lder continuous, and thus uniformly continuous on $[0,\tilde{T}]$, thus choosing $\tilde{\delta}$ sufficiently small, we can make $|Y(u)-Y(\lambda(u))|$ arbitrarily small, and therefore, as in equation (\ref{firstLambdaBound}), we can make $\int_0^s|\is_{Y(u)\geq 0}-\is_{Y(\lambda(u))\geq 0}|du$ arbitrarily small. It then follows by the triangle inequality that for any $s\in [0,\tilde{T}]$
    \begin{align} \label{secondLambdaBound} \left| \int_0^s \is_{X(u)\geq 0}du - \int_0^s \is_{Y(u)\geq 0}du \right| < \epsilon. \end{align}
    
    Now we define a function $\tilde{\lambda}:[0,T]\to \B{R}$ as follows:
    $$s(t):=\inf\left\{u\geq 0:\int_0^u\is_{X(v)\geq 0}dv\geq t\right\}\quad \text{and}\quad \tilde{\lambda}(t):=\int_0^{\lambda(s(t))}\is_{Y(u)\geq 0}du.$$

    Since $\int_0^{\tilde{T}}\is_{X(v)\geq 0}du > T$, we have that $s:[0,T]\to [0,\tilde{T}]$, and so the function $\tilde{\lambda}$ is well defined on $[0,T]$. We can also observe that $\tilde{\lambda}$ is by definition non-decreasing. Now by the triangle inequality
    \begin{align}
    |t-\tilde{\lambda}(t)| &= \left| t - \int_0^{\lambda(s(t))}\is_{Y\geq 0}\right| = \left| \int_0^{s(t)} \is_{X\geq 0} - \int_0^{s(t)}\is_{Y\geq 0} + \int_0^{s(t)}\is_{Y\geq 0} - \int_0^{\lambda(s(t))}\is_{Y\geq 0} \right| \notag \\
    & \label{thirdLineSko}\leq \left|\int_0^{s(t)}\is_{X\geq 0}-\int_0^{s(t)}\is_{Y\geq 0}\right| + \left|\int_0^{s(t)}\is_{Y\geq 0}-\int_0^{\lambda(s(t))}\is_{Y\geq 0}\right|
    \end{align}
    where we use the fact that by definition $\int_0^{s(t)} \is_{X \geq 0} = t$. Then first term of (\ref{thirdLineSko}) can be made arbitrarily small by (\ref{secondLambdaBound}) and the second term is bounded by $\tilde{\delta}$, therefore we can choose $\tilde{\delta}$ sufficiently small so that $\sup_{t\in [0,T]}|t-\tilde{\lambda}(t)|<\epsilon$.

    Next we prove that for all $t\in [0,T]$, we can make $|g(X)(t)-g(Y)(\tilde{\lambda}(t))|$ arbitrarily small by choosing $\delta$ sufficiently small. So fix $t\in [0,T]$, and define 
    $$v(t):=\inf\left\{s\geq 0: \int_0^s \is_{Y\geq 0} \geq \int_0^{\lambda(s(t))} \is_{Y\geq 0}\right\}.$$
    
    Observe that by definition of $g,\lambda, \tilde{\lambda}$ we have $g(X)(t)=X(s(t))$ and $g(Y)(\tilde{\lambda}(t))=Y(v(t))$. We can note also that by definition, $X(s(t))\geq 0$, and if $Y(\lambda(s(t)))\geq 0$, then $v(t)=\lambda(s(t))$. If $Y(\lambda(s(t)))< 0$, then by (\ref{firstLambdaBound}), $|v(t)-\lambda(s(t))|<\epsilon$. Therefore as $Y$ is, by assumption, $1/4$-H\"{o}lder continuous, and by definition of $\lambda$, $|X(s(t))-Y(\lambda(s(t))|<\epsilon$, we have that:
    \begin{align} \label{gSmall}|g(X)(t) - g(Y)(\tilde{\lambda}(t))| \leq |X(s(t))-Y(\lambda(s(t))| + |Y(\lambda(s(t))-Y(v(t))|\leq \epsilon + C|\lambda(s(t))-v(t)|^{1/4} < \epsilon + C\epsilon^{1/4},\end{align}
    which can be made arbitrarily small. 
    
    Hence we can see that $\tilde{\lambda}$ is such that $\sup_{t\in [0,T]}|t-\tilde{\lambda}(t)|$ and $\sup_{t\in [0,T]}|X(t)-Y(\tilde{\lambda}(t))|$ can be made arbitrarily small. However for Skorokhod continuity, we also require $\tilde{\lambda}$ to be increasing and bijective $[0,T]\to [0,T]$, so it remains to show that we can approximate $\tilde{\lambda}$ by an increasing bijection $[0,T]\to [0,T]$. 
    
    Since $\tilde{\lambda}:[0,T]\to \B{R}$ is increasing, there are at most countably many discontinuities. Call the points of discontinuity $t_1,t_2,\ldots$ and suppose they have jumps of size $J_i:=\tilde{\lambda}(t_i)-\tilde{\lambda}(t_i-)$. Let $N$ be such that $\sum_{i=1}^{N-1}J_i>1-\epsilon$. Then define $\hat{\lambda}_{N,n}$, $n\geq 0$, to be the function $\tilde{\lambda}$ with the discontinuities at $t_{N+1},t_{N+2},\ldots t_{N+n}$ removed:
    $$ \hat{\lambda}_{N,n}(t):= \tilde{\lambda}(t) - \sum_{i=1}^n J_{N+i}\is_{t\geq t_{N+i}} \leq \hat{\lambda}_{N,n-1}.$$
    
    Then since $(\hat{\lambda}_{N,n})_{n\in \B{N}}$ is a decreasing sequence of non-decreasing functions which is bounded below by $0$, thus $\hat{\lambda}_N:=\lim_{n\to\infty} \hat{\lambda}_{N,n}$ exists and is non-decreasing with $N$ jump discontinuities at $t_1,t_2,\ldots t_N$. Furthermore, by construction we have that $\sup_{t\in [0,T]} |\tilde{\lambda}(t)-\hat{\lambda}_N(t)|\leq \epsilon$.

    Now consider $\hat{\delta}>0$ taken sufficiently small so that the intervals $[t_i-\hat{\delta},t_i + \hat{\delta}]$, $i=1,\ldots, N$ are pairwise disjoint. Then define $\hat{\lambda}(t)$ to be the function which takes value $\hat{\lambda}_N(t)$ on $[0,T]\setminus\bigcup_{i=1}^N [t_i-\hat{\delta},t_i+\hat{\delta}]$ and which linearly interpolates across each interval. Therefore $\hat{\lambda}$ is continuous, and choosing $\hat{\delta}$ sufficiently small, we can make $\sup_{t\in [0,T]}|\tilde{\lambda}(t)-\hat{\lambda}(t)|$ arbitrarily small.

    Finally, define $\lambda_{\text{cand}}:=(\hat{\lambda}(t)+\epsilon t)/(\hat{\lambda}(T)+T\epsilon)$. So $\lambda_{\text{cand}}$ is continuous and strictly increasing bijection $\lambda_{\text{cand}}:[0,T]\to [0,T]$, and we can choose $\delta$ sufficiently small so that $\sup_{0\leq t\leq T}|\lambda_{\text{cand}}(t)-\tilde{\lambda}(t)|<\epsilon$. Therefore by the triangle inequality, we can make $\sup_{t\in [0,T]}|\lambda_{\text{cand}}(t)-t|$ arbitrarily small. Moreover, since $g(Y)(t)$ is $1/4$-H\"{o}lder continuous as a function of $t$, thus by the triangle inequality and (\ref{gSmall}):
    $$\sup_{0\leq t\leq T}|g(X)(t)-g(Y)(\lambda_{\text{cand}}(t))|\leq \sup_{0\leq t\leq T}|g(X)(t)-g(Y)(\tilde{\lambda}(t))| + \sup_{0\leq t\leq T}|g(Y)(\tilde{\lambda}(t))-g(Y)(\lambda_{\text{cand}}(t))|,$$
    can be made arbitrarily small. Therefore for $\delta$ sufficiently small, we have that $d_\infty(X,Y)<\delta$ implies $d_T(g(X),g(Y))<\epsilon$. 

    Specifically, say that for each $T\in \B{N}$, $d_\infty(X,Y)<\delta_T \implies d_T(g(X),g(Y))<\epsilon$. Thus for $d_\infty(X,Y)<\min_{T\leq N}\delta_T$, we have:
    $$
    d_\infty(g(X),g(Y))=\sum_{T=1}^\infty 2^{-T}(1\wedge d_T(g(X),g(Y)) < \sum_{T=1}^N 2^{-T}\epsilon + \sum_{T=N+1}^\infty 1 < \epsilon + 2^{-N}.
    $$
    Since $\epsilon,N$ are arbitrary, this can be made arbitrarily small. 

    Hence we can finally conclude that, for all $\epsilon > 0$, there exists $\delta>0$ such that $d_\infty(X,Y)<\delta$ implies that $d_\infty(g(X),g(Y))<\epsilon$. Therefore $g$ is continuous at $Y$. Since Brownian motion almost surely has the properties we asked of $Y$, Brownian motion $B$ is almost surely not in the discontinuity set of the transformation $g$; that is, $\B{P}(B\in D_g)=0$, as required.
\end{proof}

Recall Proposition \ref{asympSmallTC}:

\begin{propn-non} Let $(W(t))_{t\geq 0}$ be a stochastic process such that $(m^{-1/2}W(mt))_{t\geq 0}$ converges in the distribution as $m\to\infty$ to $(B(t))_{t\geq 0}$ (resp. $(|B(t)|)_{t\geq 0}$) in the Skorohod topology on $\C{D}([0,\infty),\B{R})$, where $B$ is a standard  Brownian motion. Suppose that for every $T>0$ we have
$$\sup_{0\leq t\leq T}|m^{-1}\alpha(mt)|\xrightarrow[m\to\infty]{\B{P}} 0.$$
Then $\left(m^{-1/2}W(mt+\alpha(mt))\right)_{t\geq 0}\xrightarrow[m\to\infty]{d}(B(t))_{t\geq 0}$ (resp. $(|B(t)|)_{t\geq 0}$) in the Skorohod topology on $\C{D}([0,\infty),\B{R})$
\end{propn-non}

\begin{proof} We begin by showing that the result holds for convergence in $\C{D}([0,T],\B{R})$. That is, fix $T>0$ and assume that $(m^{-1/2}X(mt))_{0\leq t\leq T} \xrightarrow[m\to\infty]{d}(B(t))_{0\leq t\leq T}$ in the Skorohod topology on $\C{D}([0,T],\B{R})$. In order to prove that $\left(m^{-1/2}X(mt+\alpha(mt))\right)_{0\leq t\leq T}\xrightarrow[m\to\infty]{d}(B(t))_{0\leq t\leq T}$, it suffices to prove that 
$$(U_m(t))_{0\leq t\leq T}:=\left( m^{-1/2}X(mt+\alpha(mt))-m^{-1/2}X(mt)\right)_{0\leq t\leq T} \xrightarrow[m\to\infty]{\B{P}} 0$$
in the Skorohod topology on $\C{D}([0,T],\B{R})$, since by Slutsky's theorem, if $A_n\xrightarrow{d} A$ and $B_n\xrightarrow{\B{P}}0$, then $A_n + B_n\xrightarrow{d}A$. Therefore we want to prove that $\forall \; \varepsilon > 0$, there exists $M$ such that $m>M$ implies 
\begin{equation}\label{aimConv}\B{P}\left(\sup_{0\leq t\leq T}|U_m(t)|>\varepsilon\right)< \varepsilon.\end{equation}

Now define the modulus of continuity of a process $Y(t)$ as $w(Y,\delta):=\underset{|s-t|\leq \delta}{\sup}|Y(s)-Y(t)|$. Note that $w(Y,\delta)=\underset{|s-t|\leq \delta}{\sup}|Y(s)-Z(s)+Z(s)-Z(t)+Z(t)-Y(t)|$ so using the triangle inequality, we have that $w(Y,\delta)\leq 2\sup_{0\leq t\leq 1}|Y(t)-Z(t)| + w(Z,\delta)$. Therefore by symmetry, it follows that $|w(Y,\delta)-w(Z,\delta)|\leq 2\sup_{0\leq t\leq 1}|Y(t)-Z(t)|$. Hence $w(Y,\delta)$ is continuous in $Y$ with respect to the Skorohod metric. Therefore:
$$\lim_{m\to\infty}w\left((m^{-1/2}X(mt))_{0\leq t\leq T},\delta\right)=w\left(\lim_{m\to\infty}(m^{-1/2}X(mt))_{0\leq t\leq T},\delta\right)\overset{d}{=}w((B(t))_{0\leq t\leq T},\delta).$$

The Portmanteau theorem (see Lemma 2.2 of \cite{asympStats}) and Levy's modulus of continuity theorem then imply that 
$$\underset{m\to\infty}\limsup \,\B{P}\left(w\left((m^{-1/2}X(mt))_{0\leq t\leq T},\delta\right)\geq \varepsilon\right)\leq \B{P}\left(w\left((B(t))_{0\leq t\leq T},\delta\right)\geq \varepsilon\right)\xrightarrow[\delta\to 0]{}0.$$

Now consider the event $E(m,\delta)$ on which $\sup_{0\leq t\leq T}|m^{-1}\alpha(mt)|>\delta$. By definition of modulus of continuity, we have that $\sup_{0\leq t\leq 1}|U_m(t)|\leq w\left((m^{-1/2}X(mt))_{0\leq t\leq T},\delta\right)$ on $E(m,\delta)^c$, and therefore:
$$
\limsup_{m\to\infty} \, \B{P}\left( \sup_{0\leq t\leq T} |U_m(t)| \mathbbm{1}_{E(m,\delta)^c}>\varepsilon\right)\xrightarrow[\delta\to 0]{}0
$$

Alternatively, since $\sup_{0\leq t\leq T}|m^{-1}\alpha(mt)|\xrightarrow[m\to\infty]{\B{P}}0$, we have that $\B{P}(E(m,\delta))\to 0$ as $m\to \infty$, and hence $\B{P}\left(\sup_{0\leq t\leq T} |U_m(t)|\mathbbm{1}_{E(m,\delta)}>\varepsilon\right)\leq \B{P}\left(E(m,\delta) \right) \xrightarrow[m\to\infty]{}0$. Therefore
$$
\lim_{m\to\infty}\B{P}\left( \sup_{0\leq t\leq T} |U_m(t)|\mathbbm{1}_{E(m\delta)}>\varepsilon\right) = 0
$$

So  combining the convergence on $E(m,\delta)$ and $E(m,\delta)^c$, equation (\ref{aimConv}) is proven as required, and hence by Slutsky's theorem, we have 
$$\left(m^{-1/2}X(mt+\alpha(mt))\right)_{0\leq t\leq T}\xrightarrow[m\to\infty]{d}(B(t))_{0\leq t\leq T}.$$

Finally, since the convergence holds for all $T$, thus by Lemma 3 of section 16 of \cite{bills}, it holds that
$$\left(m^{-1/2}X(mt+\alpha(mt))\right)_{t\geq 0}\xrightarrow[m\to\infty]{d}(B(t))_{t\geq 0}$$
in the Skorokhod topology on $\C{D}([0,\infty),\B{R})$.
\end{proof}

\newpage

\section*{Appendix: The Asymptotic Velocity of the $N$-BBM}
The aim of this section is to prove the following theorem about the velocity of the $N$-BBM system.

\begin{theorem} \label{asympVeloc}
    $$v_N = \sqrt{2} - \frac{\pi^2}{\sqrt{2}(\log N)^2} + o\left( \frac{1}{(\log N)^2}\right).$$
\end{theorem}
Our proof will use Theorem 7.6 of \cite{barakErdos} of Mallein and Ramassamy, which we state here for completeness. The theorem gives the asymptotic speed of discrete-time branching random walks with selection ($N$-BRW). In the $N$-BRW, we have $N$ particles and at each discrete time, every particle branches ito a random number of offspring, located at random distances from their parents, and we simultaneously delete all but the $N$ rightmost particles. The offspring is described by a distribution $\C{M}$, which is a random point process. 

\begin{theorem} \textit{(Theorem 7.6 in \cite{barakErdos})} Let $Z^N$ be an $N$-BRW and let $M\sim \C{M}$ satisfy:
\begin{itemize}
    \item $\kappa(\theta):=\log\B{E}\left(\sum_{m\in M}e^{\theta m}\right) < \infty \quad \forall \theta >0 $
    \item $\exists \theta^* > 0$ such that $\theta^* \kappa'(\theta^*)=\kappa(\theta^*)$
    \item $\B{E}\left(|\max_{m\in M}m|^2\right)<\infty$
    \item $\B{E}\left(\sum_{m\in M}e^{\theta^* m}\left(\log \sum_{m\in M}e^{\theta^* m}\right)^2\right)<\infty$
\end{itemize}
Then $v_N:=\lim_{t\to \infty} t^{-1}Z^N_1(t)=\lim_{t\to\infty}t^{-1}Z^N_N(t)$ exists and $v_N=\kappa'(\theta^*) + \frac{\pi^2 \theta^* \kappa''(\theta^*)}{2(\log N)^2} + o\left((\log N)^{-2}\right)$
\end{theorem}

We can thus couple $N$-BBM to appropriately chosen $N$-BRWs to determine the speed of $N$-BBM.

\begin{lemma} \label{speedUB} $v_N \leq \sqrt{2} - \frac{\pi^2}{\sqrt{2}(\log N)^2} + o\left( \frac{1}{(\log N)^2}\right)$
\end{lemma}

\begin{proof} Define $Z^{N,1}(t)$ to be a branching Brownian motion process in which at each integer time we delete all but the $N$ rightmost particles. Therefore $Z^{N,1}$ may have more than $N$ particles at some time $t\in \B{R}\setminus \B{Z}$, but has exactly $N$ particles at each time in $\B{Z}$. 

By Proposition 3 of \cite{hydroNBBM}, there is a natural coupling between $Z^N$ and $Z^{N,1}$ such that $Z^N(0)\preccurlyeq Z^{N,1}(0) \implies Z^N(t)\preccurlyeq Z^{N,1}(t)$ for all subsequent times $t>0$. Therefore, defining $v_{N,1}:=\lim_k k^{-1}Z_1^{N,1}(k)$, then certainly $v_{N}\leq v_{N,1}$. 

By construction, taking $Z^{N,1}$ at integer times, $(Z^{N,1}(n))_{n\in \B{N}}$ is an $N$-BRW with reproduction law $\C{M}_1$, say. The law $\C{M}_1$ describes the distribution of particles at time $1$ of a branching Brownian motion process started from a single particle at the origin. Let $M_1$ be a random point process of law  $\C{M}_1$. It is easy to check that $\C{M}_1$ satisfies the condition of the above Theorem, and we can calculate $\kappa(\theta) = 1+\theta^2/2$. Therefore the conclusion of the Theorem gives $\theta^*=\sqrt{2}$ and
$$
v_N \leq v_{N,1} = \sqrt{2} - \frac{\pi^2}{\sqrt{2}(\log N)^2} + o\left( \frac{1}{(\log N)^2}\right),
$$
as required. 
\end{proof}

For the lower bound, we will define another branching random walk $(Y^{N,\delta}(t))_{t\geq 0}$ with selection at discrete times $(k\delta)_{k\in\B{N}}$. Consider a branching Brownian motion process $\hat{Y}^1$ starting with a single particle at the origin in which we permit at most one branching event. Therefore there at most 2 particles in $\hat{Y}^1$ at any time. Let the distribution of the particles of $\hat{Y}^1(\delta)$ be called $\C{\hat{M}}_\delta$. Now define $(Y^{N,\delta}(n\delta))_{n\in\B{N}}$ to be the $N$-BRW starting with $\lfloor \frac{N}{2} \rfloor$ particles with offspring distribution $\C{\hat{M}}_\delta$ and with selection of the right-most $\lfloor \frac{N}{2} \rfloor$ particles. Then since there are always fewer particles in the process $Y^{N,\delta}$ than in $Z^{N}$, Lemma 7.4 of \cite{barakErdos} proves that there exists a coupling such that $Y^{N,\delta}(0)\oop Z^{N}(0) \implies Y^{N,\delta}(t)\oop Z^N(t)$ for all $t\geq 0$. Denoting the speed of the process $Y^{N,\delta}$ by $\hat{v}_{N,\delta}$, it immediately follows that $\hat{v}_{N,\delta}\leq v_N$ Thus it only remains give the velocity $\hat{v}_{N,\delta}$.

\begin{lemma} $v_N \geq \sqrt{2} - \frac{\pi^2}{\sqrt{2}(\log N)^2} + o\left( \frac{1}{(\log N)^2}\right)$
\end{lemma}

\begin{proof} Let $\hat{M}_\delta$ be a random point process with law $\C{\hat{M}}_\delta$. Conditioning on the branching time, we calculate:
\begin{align*}
    \hat{\kappa}_\delta(\theta)&:=\log\B{E}\left(\sum_{m\in \hat{M}_\delta}e^{\theta m} \right) = \log\left( \int_0^\infty \B{E}\left[\sum_{m\in \hat{M}_\delta}e^{\theta m}\middle| \tau=t\right]e^{-t}dt\right) \\
    &=\log \left( \int_0^\delta 2e^{-t}\B{E}[e^{\theta B_\delta}]dt + \int_\delta^\infty e^{-t}\B{E}[e^{\theta B_\delta}]dt\right) \\
    &=\frac{\theta^2 \delta}{2} + \log (2-e^{-\delta})
\end{align*}
where $(B_t)_{t\geq 0}$ is a standard Brownian motion, and therefore $\theta^* = \sqrt{2\delta^{-1}\log(2-e^{-\delta})}$. So the above Theorem gives
$$
v_N \geq \hat{v}_{N,\delta} = \sqrt{2\delta^{-1}\log(2-e^{-\delta})} - \frac{\pi^2 \sqrt{2\delta^{-1}\log(2-e^{-\delta})}}{\sqrt{2}(\log N)^2} + o\left(\frac{1}{(\log N)^2}\right)
$$
and by L'H\^{o}pital's rule, $\sqrt{2\delta^{-1}\log(2-e^{-\delta})} \to \sqrt{2}$ as $\delta \downarrow 0$, so $v_N \geq \sqrt{2}-\frac{\pi^2}{2(\log N)^2} + o((\log N)^{-2})$ as required. 
\end{proof}

\begin{proof} \textit{(of Theorem \ref{asympVeloc})}
Putting together the preceding two lemmas gives the desired result immediately. 
\end{proof}


\section*{Acknowledgements}
The author thanks Julien Berestycki for his supervision of the project and heplful feedback on the paper, as well as helpful conversations with Bastien Mallein, and constructive feedback from Matthias Winkel and Brett Kolesnik.

\bibliographystyle{plain}
\bibliography{bps}

\subfile{}

\end{document}